\documentclass[11pt,a4paper]{article}
\usepackage{graphicx,amsfonts,amsbsy, amssymb, amsthm, amsmath}
\usepackage{bbm}
\usepackage{fouridx}
\usepackage{pgfplots}

\usepackage{authblk}

\renewcommand{\vec}[1]{{\mathbf{#1}}}

\renewcommand{\Re}{\mathop{\mathfrak{Re}}}

\newcommand{\rmd}{{\mathrm d}}
\newcommand{\rme}{{\mathrm e}}
\newcommand{\rmi}{{\mathrm i}}

\newcommand{\Ord}{{\mathrm O}}
\newcommand{\littleo}{{\mathrm o}}

\DeclareSymbolFont{lettersA}{U}{pxmia}{m}{it}

\DeclareMathSymbol{\alphaup}{\mathord}{lettersA}{"0B}
\DeclareMathSymbol{\betaup}{\mathord}{lettersA}{"0C}
\DeclareMathSymbol{\gammaup}{\mathord}{lettersA}{"0D}
\DeclareMathSymbol{\deltaup}{\mathord}{lettersA}{"0E}
\DeclareMathSymbol{\epsilonup}{\mathord}{lettersA}{"22}
\DeclareMathSymbol{\zetaup}{\mathord}{lettersA}{"10}
\DeclareMathSymbol{\etaup}{\mathord}{lettersA}{"11}
\DeclareMathSymbol{\thetaup}{\mathord}{lettersA}{"12}
\DeclareMathSymbol{\iotaup}{\mathord}{lettersA}{"13}
\DeclareMathSymbol{\kappaup}{\mathord}{lettersA}{"14}
\DeclareMathSymbol{\lambdaup}{\mathord}{lettersA}{"15}
\DeclareMathSymbol{\muup}{\mathord}{lettersA}{"16}
\DeclareMathSymbol{\nuup}{\mathord}{lettersA}{"17}
\DeclareMathSymbol{\xiup}{\mathord}{lettersA}{"18}
\DeclareMathSymbol{\piup}{\mathord}{lettersA}{"19}
\DeclareMathSymbol{\rhoup}{\mathord}{lettersA}{"1A}
\DeclareMathSymbol{\sigmaup}{\mathord}{lettersA}{"1B}
\DeclareMathSymbol{\tauup}{\mathord}{lettersA}{"1C}
\DeclareMathSymbol{\upsilonup}{\mathord}{lettersA}{"1D}
\DeclareMathSymbol{\phiup}{\mathord}{lettersA}{"1E}
\DeclareMathSymbol{\chiup}{\mathord}{lettersA}{"1F}
\DeclareMathSymbol{\psiup}{\mathord}{lettersA}{"20}
\DeclareMathSymbol{\omegaup}{\mathord}{lettersA}{"21}

\newcommand{\vecphi}{{\pmb{\phiup}}}

\renewcommand{\Psi}{\varPsi}
\renewcommand{\Lambda}{\varLambda}
\renewcommand{\Sigma}{\varSigma}
\renewcommand{\Gamma}{\varGamma}
\renewcommand{\Theta}{\varTheta}
\renewcommand{\Xi}{\varXi}
\renewcommand{\Pi}{\varPi}
\renewcommand{\Upsilon}{\varUpsilon}
\renewcommand{\Phi}{\varPhi}
\renewcommand{\Omega}{\varOmega}

\newcommand{\R}{{\mathbb R}}
\newcommand{\N}{{\mathbb N}}

\newcommand{\Z}{{\mathbb Z}}

\newcommand{\C}{{\mathbb C}}

\newcommand{\Prob}{{\mathbb P}}

\newcommand{\coloneq}{\mathbin{\hbox{\raise0.08ex\hbox{\rm :}}\!\!=}}
\newcommand{\eqcolon}{\mathbin{=\!\!\hbox{\raise0.08ex\hbox{\rm :}}}}

\renewcommand{\leq}{\leqslant}
\renewcommand{\geq}{\geqslant}
\renewcommand{\epsilon}{\varepsilon}

\newcommand{\dimostrazione}{\noindent{\sl Proof.}\phantom{X}}
\newcommand{\dimostrazionea}[1]{\noindent{\sl Proof of #1.}\phantom{X}}
\newcommand{\finire}{\hspace*{\fill}~$\Box$}

\newcommand \printdate[3]{%
    \def \@suffix##1{%
        \def \@n{##1}%
        \ifnum \@n = 1 st\else%
        \ifnum \@n = 2 nd\else%
        \ifnum \@n = 3 rd\else%
        \ifnum \@n = 21 st\else%
        \ifnum \@n = 22 nd\else%
        \ifnum \@n = 23 rd\else%
        \ifnum \@n = 31 st\else%
        th\fi \fi \fi \fi \fi \fi \fi%
    }%
    \relax%
    \number #1\raise0.7ex\hbox{\footnotesize \@suffix{#1}}\kern0.25em%
    \ifcase #2\or%
        January\or February\or March\or%
        April\or May\or June\or%
        July\or August\or September\or%
        October\or November\or December%
    \fi\ %
    \number #3%
}

\newtheorem{theorem}{Theorem}[section]
\newtheorem{proposition}[theorem]{Proposition}

\newtheorem{corollary}[theorem]{Corollary}
\newtheorem{lemma}[theorem]{Lemma}

\newcommand{\rmD}{\mathrm{D}}
\newcommand{\rmE}{\mathrm{E}}

\newcommand{\twoFone}[1]{\fourIdx{}{2}{(#1)}{1}{F}}
\newcommand{\1}{\vec{1}}
\newcommand{\JBE}{J$\betaup$E}

\numberwithin{equation}{section}

\usepackage[text={15cm,23cm}]{geometry} 

\usepackage[colorlinks=true,%
            linkcolor=red!50!black,%
            citecolor=blue!50!black,%
            urlcolor=darkgray]{hyperref}  

\begin{document}
\title{Extreme eigenvalues of random matrices from Jacobi ensembles}
\author{B.~Winn}
\affil{Department of Mathematical Sciences, 
Loughborough University, Loughborough,
LE11 3TU, U.K.}
\date{\printdate{22}1{2024}}
\maketitle
\begin{abstract}
Two-term asymptotic formul\ae\ for the probability distribution functions
for the smallest eigenvalue of the Jacobi $\beta$-Ensembles are derived 
for matrices of large size in the r\'egime where $\beta>0$ is arbitrary
and one of the model parameters $\alpha_1$ is an integer.  By a straightforward
transformation this leads to corresponding results for the
distribution of the largest eigenvalue.  The explicit expressions are
given in terms of multi-variable hypergeometric functions, and it is found that
the first-order corrections are proportional to the derivative of the
leading order limiting distribution function.

In some special cases $\beta=2$ and/or small values of $\alpha_1$, 
explicit formul\ae\ involving more familiar functions, such as the 
modified Bessel function of the first kind, are presented.
\end{abstract}

\thispagestyle{empty}

\section{Introduction}
A random matrix is a matrix whose entries are random variables.  As
eigenvalues of a matrix are continuous functions of its entries, so
the eigenvalues of a random matrix are random variables.
A random $N\times N$ matrix has the Jacobi $\beta$-Ensemble (\JBE) distribution
if a joint probability density function of its eigenvalues is
\begin{equation}
\label{eq:2}
\frac1{S_N(\alpha_1+1, \alpha_2+1, \beta/2)}  \prod_{i=1}^N
x_i^{\alpha_1} (1-x_i)^{\alpha_2}|\Delta(\vec{x})|^{\beta},\qquad
0\leq x_i \leq 1,
\end{equation}
where
\begin{equation}
\label{eq:16}
  S_N(a ,b ,c) \coloneq \prod_{i=0}^{N-1} \frac{\Gamma(1+(i+1)
c) \Gamma(a +ic)\Gamma(b +ic)}{\Gamma(1+c)\Gamma(
a +b +(N+i-1)c)},
\end{equation}
and the Vandermonde determinant is defined by
\begin{equation}
  \label{eq:1}
  \Delta(\vec{x}) \coloneq \prod_{1\leq i < j \leq N} (x_j-x_i).
\end{equation}
That \eqref{eq:2} is a properly normalised probability density
is a consequence of Selberg's integral \cite{sel:boe}.

In many situations $\beta$ is a non-negative integer, and we
will mostly be assuming that one of $\alpha_1$ or $\alpha_2$ is
a non-negative integer, but \eqref{eq:2} makes sense for
arbitrary real values of these parameters subject to the constraints
\begin{equation}
  \label{eq:3}
  \alpha_1 > -1,\quad \alpha_2 > -1, \quad \beta > - \frac12\min\left\{
\frac1N, \frac{\alpha_1}{N-1}, \frac{\alpha_2}{N-1}\right\}.
\end{equation}
The naming of this ensemble reflects the presence in
\eqref{eq:2} of the factors $x_i^{\alpha_1}(1-x_i)^{\alpha_2}$
which are a density with respect to which (a certain version of) the 
classical Jacobi polynomials form an orthogonal collection.

We label by $\phi_i$ the sorted eigenvalues, so that
$0 \leq \phi_1 \leq \phi_2 \leq \cdots \leq \phi_N \leq 1$.  This
article is concerned with the distribution of the extreme
eigenvalues $\phi_1$ and $\phi_N$.  In fact, since the change of
variables $x_i \mapsto 1 - x_i$, for $i=1,\ldots,N$, in \eqref{eq:2} 
leaves the joint probability density invariant, save for the exchange
$\alpha_1 \leftrightarrow \alpha_2$, and reverses the order of the eigenvalues,
it will not present a loss of generality to focus on the 
\emph{smallest} eigenvalue $\phi_1$.

The limiting empirical eigenvalue density for Jacobi random matrices
was derived in \cite{wac:tle}. For fixed $\alpha_1,\alpha_2$, the large $N$ 
limiting density is
\begin{equation}
  \label{eq:151}
  \frac{N}\pi \frac1{\sqrt{x(1-x)}},\qquad 0<x<1.
\end{equation}
This means that for large $N$ the number of eigenvalues in the
interval  $[0,N^{-2}]$ is approximately
$N/\pi\int_0^{N^{-2}}(x(1-x))^{-1/2}\,\rmd x = \Ord(1)$,
and it is natural to expect $N^2\phi_1$ to have a non-trivial limiting
distribution.

Our main objects of interest will be the (cumulative) probabilty 
distribution function $F_{\phi_1}(\xi) \coloneq \Prob(\phi_1 \leq \xi)$,
and the rescaled version
\begin{equation}
  \label{eq:8}
F_{N^2\phi_1}(x) = \Prob( N^2 \phi_1 \leq x) = 
\Prob\left( \phi_1 \leq \frac{x}{N^2}\right) = 
F_{\phi_1}\left( \frac{x}{N^2} \right).
\end{equation}
Deferring to below a more comprehensive summary of previous work on
this problem, we mention a result \cite{mor:eed} of 
 Moreno-Pozas, Morales-Jimenez, McKay in the case $\beta=2$ (the
Jacobi Unitary Ensemble, JUE).  They proved, for
$\alpha_1=0,1$ and $\alpha_2\in\N_0$; and for $\alpha_1=2,
\alpha_2\in\{0,1,2\}$, the two-term asymptotic result
\begin{equation}
\label{eq:9}
  F_{N^2\phi_1}(x) = 1 - \rme^{-x}\det(I_{j-i}(2\sqrt{x})) + 
\frac{\alpha_1+\alpha_2}N x\rme^{-x}
\det(I_{2+j-i}(2\sqrt x)) + \Ord\left(\frac1{N^2}\right)
\end{equation}
where the determinants appearing in~\eqref{eq:9} are of size 
$\alpha_1\times\alpha_1$, and $I_n(z)$ is the $I$-Bessel
function---the modified Bessel function of the first kind,
\begin{equation}
  \label{eq:10}
I_\nu(z) \coloneq \frac{z^\nu}{2^\nu\sqrt{\pi}\Gamma(\nu+\frac12)}\int_{-1}^1
\rme^{-zt}(1-t^2)^{\nu-1/2}\,\rmd t,\qquad\Re\{\nu\}>-\frac12,
\end{equation}
and $I_{-n}(z)=(-1)^n I_n(z)$ for $n\in\Z$.

On the other hand, Borodin and Forrester have derived \cite{bor:isa}
the leading-order distribution of the smallest eigenvalue of the \JBE\
for any $\beta>0$ and $\alpha_1\in\N_0$:
\begin{equation}
\label{eq:150}
  \lim_{N\to\infty} F_{N^2\phi_1}(x)
=
1 - \rme^{-\beta x/2} 
 \fourIdx{}0{(\beta/2)}1{F}\left(;\frac{2\alpha_1}{\beta};x\vec{1}^{\alpha_1}
\right),
\end{equation}
where $\fourIdx{}{0}{(\sigma)}{1}{F}(;c;\vec{x})$ is a multivariate 
hypergeometric function that will be defined precisely 
in Section~\ref{sec:mult-hyperg-funct}, and 
$\1^n\coloneq(1,1,\ldots,1)\in\R^n$.  Our principal result is a
version of the two-term asymptotic \eqref{eq:9} valid for $\beta>0$.
\begin{theorem}  \label{thm:main}
Let $\phi_1$ be the smallest eigenvalue of the $N\times N$ Jacobi
$\beta$-Ensemble, $\beta>0$, with
$\alpha_1\in\N_0$ and $\alpha_2>-1$.  For $x>0$,
  \begin{multline}
  \label{eq:140}%
F_{N^2\phi_1}(x) = 1 - \rme^{-\beta x/2} \fourIdx{}0{(\beta/2)}1F\left(
;\frac{2\alpha_1}\beta; x\vec{1}^{\alpha_1} \right) \\
 + \frac{x^{1+\alpha_1}}N\left((\alpha_1 + \alpha_2 +1)-\frac\beta2\right) 
\left(\frac\beta2\right)^{2\alpha_1} 
\frac{\Gamma(1+\beta/2)}{\Gamma(1+\alpha_1)\Gamma(1+\alpha_1+\beta/2)}
\\\times \rme^{-\beta x/2}
\fourIdx{}0{(\beta/2)}1F\left(
;\frac{2\alpha_1}\beta+2; x\vec{1}^{\alpha_1} \right) +
\Ord\left(\frac1{N^2}\right).
\end{multline}
The error estimate can depend on $\alpha_1, \alpha_2, \beta$ but is
uniform for $x$ in a compact set.
\end{theorem}
All our results for the distribution of the smallest eigenvalue
can be re-cast to give an analogous result for the largest
eigenvalue, as indicated earlier.  We will not write down these
analogues for every result, allowing just the following Corollary
of Theorem \ref{thm:main}.
\begin{corollary}  \label{cor:main}
Let $\phi_N$ be the largest eigenvalue of the $N\times N$ Jacobi
$\beta$-Ensemble, $\beta>0$, with
$\alpha_2\in\N_0$ and $\alpha_1>-1$.  For $x>0$,
  \begin{multline}
  \label{eq:140b}%
\Prob(\phi_N \leq 1 - x/N^2)
=  \rme^{-\beta x/2} \fourIdx{}0{(\beta/2)}1F\left(
;\frac{2\alpha_2}\beta; x\vec{1}^{\alpha_2} \right) \\
 - \frac{x^{1+\alpha_2}}N\left((\alpha_1 + \alpha_2 +1)-\frac\beta2\right) 
\left(\frac\beta2\right)^{2\alpha_2} 
\frac{\Gamma(1+\beta/2)}{\Gamma(1+\alpha_2)\Gamma(1+\alpha_2+\beta/2)}
\\\times \rme^{-\beta x/2}
\fourIdx{}0{(\beta/2)}1F\left(
;\frac{2\alpha_2}\beta+2; x\vec{1}^{\alpha_2} \right) +
\Ord\left(\frac1{N^2}\right).
\end{multline}
\end{corollary}

There are several known random matrix models that lead to \JBE\ 
eigenvalue distributions.  Most famous are perhaps the double-Wishart
(or \textsc{Manova}) models from Statistics:
 set $M_1$, $M_2$ to be independent $n_1\times N$
and $n_2\times N$ matrices with independent standard normal real random
variable entries, $n_1, n_2\geq N$. If
 $A=M_1^\dag M_1$, $B= M_2^\dag M_2$, then the
matrix $A(A+B)^{-1}$ has eigenvalues distributed according to the 
\JBE\ with $\beta=1$, $\alpha_1=(n_1-N-1)/2$ and $\alpha_2=(n_2-N-1)/2$
\cite{fis:tsd,hsu:otd,roy:psa,gir:ots,moo:otd}.
Since our results rely on $\alpha_1$ being integer, this requires
$n_1-N$ to be an \emph{odd} diference.

If we repeat the above construction, with complex normal random
variables, then the eigenvalue distribution of $A(A+B)^{-1}$ is
\JBE\ with $\beta=2$, $\alpha_1=n_1-N$ and $\alpha_2=n_2-N$
\cite[Section 8]{jam:dom}.

Another model leading to the joint probability density function
\eqref{eq:2} is the corners process of random matrices from classical
compact groups: if $U$ is a random $m\times m$ unitary or orthogonal matrix
chosen with respect to Haar measure, $m\geq 2N$, and $M$ is the
principal $N\times N$ submatrix of $U$ (the upper-left corner matrix), 
then, letting $s_1,\ldots,s_N$
denote the $N$ eigenvalues of $M^\dag M$,
the points $x_1=s_1/\beta,\ldots,x_N=s_N/\beta$ are distributed 
according to \eqref{eq:2} with $\alpha_1=\beta/2-1$
$\alpha_2=(m-2N+1)\beta/2-1$ and $\beta=2$ (unitary case) or $\beta=1$
(orthogonal case) \cite[\S7.2]{col:por,eat:gia}.
In the latter case $\alpha_1=-1/2$ which is not
an integer, so Theorem \ref{thm:main} does not apply, but
Corollary \ref{cor:main} does apply
for the distribution of the largest eigenvalue
if $m$ is an odd number (whence $\alpha_2$ is an integer).

Random matrix models that allow full exploration of the parameter space, 
including to arbitrary $\beta>0$, are also known
\cite{lip:amm,kil:mmf,ede:tbj}.

The \JBE\ exhibits two ``hard edges'' in the spectrum at $x=0$ and $x=1$
since the eigenvalues are strictly confined between these values, which
furthermore coincides with the support of the limiting eigenvalue 
density \eqref{eq:151}.
This is in contrast to some other random matrix models such as the
Gaussian ensembles \cite{meh:rm} which have compactly supported
limiting eigenvalue density---the famous Wigner's semi-circle law
\cite{wig:cvo,wig:otd}---but without any intrisic obstacle to
having individual eigenvalues appearing at any point on the real line.
Statistics such as the distribution of smallest eigenvalues are expected
to be ``universal'' in the limit $N\to\infty$, in the sense that they ought 
not to depend on the precise features of the random matrix model in
question.  In our present context it means that the limiting distribution
\eqref{eq:150} will be valid for other matrix models with a hard 
spectral edge.  Indeed, the same limiting distribution has been proven
for a different set of matrix models exhibiting a hard edge---the 
Laguerre $\beta$-Ensembles (L$\betaup$E; sometimes called
Wishart random matrices) \cite{for:era}, as well
as modifications of the JUE that preserve the hard edge \cite{kui:ufe}. 

The finite $N$ corrections to the leading order derived in Theorem
\ref{thm:main} are not expected to be universal---indeed the presence
of the parameter $\alpha_2$ seems to rule that out---but 
they do exhibit an interesting feature that had already been conjectured 
for the Laguerre Unitary Ensemble at the hard edge \cite{ede:bui} and 
proved for that model in \cite{bor:ano,per:fNc,hac:lcc}:
the correction term is proportional to the \emph{derivative} of
the main term.  This holds for our two-term asymptotic \eqref{eq:140},
although it may not seem immediately apparent: see \eqref{eq:106} below.
Forrester and Trinh \cite{for:fsc} have investigated the eigenvalue
density for the L$\betaup$E for $\beta>0$,
and found two-term asymptotics at the hard edge of the spectrum, and
that the correction term is also proportional to the derivative of the
leading term.  It seems likely that the methods in the present work
could also be adapted to study the hard-edge of the L$\betaup$E too.

\JBE\ random matrices have a number of known applications.  The outage
probability of multiple-input/multiple-output (MIMO) systems subject
to interference, such as those used in cellular mobile radio networks,
can be modelled in terms of the largest eigenvalue of JUE matrices
\cite{kan:qfi}.  The conductance eigenvalues in random matrix models
for mesoscopic disordered quantum systems are known to be governed by
the \JBE\ distribution with $\alpha_2=0$ \cite{bee:rmt,for:qcp}. In
this context, expressions for the average spectral density have been
derived in terms of multi-variable hypergeometric functions
\cite{viv:ted}, somewhat similar to expressions for the smallest
eigenvalue derived in Section \ref{sec:main-calculations}.
Finally, some tests in multivariate Statistics are based on the
distributions of extreme eigenvalues of \JBE\ (generally the
parameters $\beta=1$ and $\beta=2$ corresponding to the real and
complex underlying fields are most relevant), see Roy \cite{roy:oah}.  
Some of these statistical applications are reviewed in
Section 2 of \cite{joh:maa}.  Roy's test has practical applications
in signal analysis in the presence of coloured noise, for which the
distribution of the largest eigenvalue of the JUE is required
\cite{cha:ebd}.

Aside from the references \cite{mor:eed,bor:isa} mentioned above,
theoretical work on the distribution of extreme eigenvalues for Jacobi
ensembles goes back at least to \cite{kha:dot} for $\beta=2$ and
Constantine \cite{con:snd} for $\beta=1$, motivated by the
aforementioned applications in Statistics.

In \cite{dum:dot} expressions were derived for distribution functions
in terms of multivariate hypergeometric functions in $N$ variables,
and corresponding formulae for density functions in $N-1$ variables
given in \cite{dum:sed} and \cite{dre:cwb}.  Algorithms for a numerical
evaluation of the distribution of the smallest eigenvalue in the
JUE were given in \cite{due:tle} with methods applicable to arbitrary
values of the parameters $\alpha_1, \alpha_2 > -1$, and furthermore
which extend even to non-integer values of $N$.

Johnstone \cite{joh:maa} and Jiang \cite{jia:ltf} have investigated
statistics of extreme eigenvalues, and other quantities, in a setting
where the parameter values $\alpha_1$ and $\alpha_2$ are not fixed,
but vary as $N\to\infty$, leading to a soft edge in the spectrum.
Scaling limits at the hard and soft-edge were treated together in
\cite{hol:eso}.

Forrester and Li \cite{for:roc} have studied eigenvalue correlations
for a broader class of unitary ensembles with a hard edge at the 
spectrum (which includes the JUE) and found $1/N$-correction terms
consistent with \cite{mor:eed}.

In Section \ref{sec:multi-vari-hyperg} we introduce some of the 
analytic tools that will be used (multi-variable hypergeometric
functions and Jacobi polynomials).  In Section \ref{sec:main-calculations}
we collect some exact formul\ae\ for finite $N$.  In Section 
\ref{sec:two-term-asymptotic} we prove a two-term asymptotic
formula for $\fourIdx{}2{(\sigma)}1{F}$
multivariate hypergeometric functions, that is then
used to give the proof of Theorem \ref{thm:main} in Section
\ref{sec:main-result}.  A few special cases are treated in
Section \ref{sec:explicit-formulas}.

\section{Multi-variable hypergeometric functions and 
Jacobi polynomials}  \label{sec:multi-vari-hyperg}
Multi-variable analogues of classical hypergeometric functions and
orthogonal polynomials are a relatively recently-developed area
of study that have nevertheless proved very useful in Random
Matrix Theory, see, e.g.\ \cite{dum:MOPS,bak:tcs,for:sds,%
win:dmf,mez:mot,jon:sft} as well as many other articles cited in the
present work.  They can be defined as series of Jack 
polynomials which we define first.
\subsection{Jack polynomials}
Let $\vec{x}=(x_1,\ldots,x_n)$ be a set of variables,
$\lambda=(\lambda_1,\ldots,\lambda_n)$ be an integer partition\footnote{%
We can assume that the number of parts of $\lambda$ is equal to the
number of variables, since if there are more parts than variables the
corresponding Jack polynomial is zero; on the other hand, any partition
can be padded with $0$s to increase the number of parts to $n$.}
of size $|\lambda|=\lambda_1+\cdots+\lambda_n$, and
let $\sigma>0$. The Jack polynomials \cite{jac:aco} are certain homogeneous,
symmetric polynomials $C_\lambda^{(\sigma)}(\vec{x})$ of degree $|\lambda|$.

We define the operators
\begin{equation}
  \label{eq:diff1}
  \rmD_k \coloneq \sum_{i=1}^n x_i^k \frac{\partial^2}{\partial x_i^2} +
  \frac2\sigma \sum_{i\neq j} \frac{x_i^k}{x_i-x_j}\frac\partial{\partial x_i}
\end{equation}
and
\begin{equation}
  \label{eq:diff2}
  \rmE_k \coloneq \sum_{i=1}^n x^k \frac\partial{\partial x_i},
\end{equation}
for $k\in\N_0$

Jack polynomials are joint eigenfunctions of $\rmE_1$ and $\rmD_2$ 
\cite{sta:scp}.  In fact,
\begin{equation}
  \label{eq:diff3}
  \rmE_1 C^{(\sigma)}_\lambda(\vec{x}) = |\lambda|C^{(\sigma)}_\lambda(\vec{x})
\end{equation}
(a relation satisfied by \emph{any} homogeneous polynomial of degree
$|\lambda|$) and
\begin{equation}
  \label{eq:diff4}
  \rmD_2 C^{(\sigma)}_\lambda(\vec{x}) = \left(\rho_\lambda +
    \frac2\sigma|\lambda|(n-1))\right)
  C^{(\sigma)}_\lambda(\vec{x}),
\end{equation}
where
\begin{equation}
  \label{eq:diff5}
  \rho_\lambda \coloneq \sum_{i=1}^n \lambda_i \left( \lambda_i - 1 -
  \frac2\sigma (i-1) \right).
\end{equation}
The definition of $C_\lambda^{(\sigma)}(\vec{x})$ is completed by 
triangularisation:
if 
\begin{equation}
  \label{eq:152}
  C_\lambda^{(\sigma)}(\vec{x}) = \sum_{\mu} b_{\mu\lambda} m_\mu(\vec{x})
\end{equation}
is the expansion of $C_\lambda^{(\sigma)}$ in the basis of
monomial symmetric functions, then the coefficient $b_{\mu\lambda}=0$
unless $\mu\leq\lambda$ in terms of dominance ordering of partitions
\cite{sta:scp};
and normalisation:
\begin{equation}
\label{eq:147}
  \sum_{|\lambda|=k} C_\lambda^{(\sigma)}(\vec{x}) = (x_1+\cdots+x_n)^k,\qquad
k\in\N.
\end{equation}
(That such a normalisation exists is proved in
\cite[Prop.~2.3]{sta:scp}, although a different normalisation for the
Jack polynomials is actually used throughout \cite{sta:scp}.  The
normalisation leading to \eqref{eq:147} is commonly-used for applications 
in Random Matrix Theory.)

\subsection{Multi-variable Jacobi polynomials}
As with multi-variable hypergeometric functions defined in the next
subsection, multi-variable generalisations of the classical Jacobi
polynomials were initially studied for Jack parameter $\sigma=2$
\cite{jam:gjp} 
with applications in Statistics in mind.  Later these were generalised
to other values of $\sigma$, with a variety of conventions for
normalisation and support of the orthogonality measure \cite{vre:ffe,
  hec:rsaI, opd:sao, ols:mjp}. They sometimes go by the name ``Jacobi
polynomials associated with the root system $BC_n$'' \cite{koo:sfa}.
In our definitions, we follow \cite{las:pdj} with a difference in the
choice of normalisation.

For $a,b\in\R$ fixed, $J_\lambda^{\sigma,a,b}(\vec{x})$ is a symmetric
polynomial eigenfunction of the operator
\begin{equation}
  \label{eq:139}
  \rmD_2 - \rmD_1 + (a+b+2)\rmE_1 - (a+1)\rmE_0
\end{equation}
of the form
\begin{equation}
  \label{eq:116}
  J_\lambda^{\sigma,a,b}(\vec{x}) = \sum_{\mu\subseteq\lambda} c_{\mu,\lambda}
C_\mu^{(\sigma)}(\vec{x}),
\end{equation}
for constants $c_{\mu,\lambda}$ depending on $a, b$ and $\sigma$, and the
notation $\mu\subseteq\lambda$ means $\mu_i\leq\lambda_i$ for $i=1,\ldots,n$.
We normalise $J_\lambda^{\sigma,a,b}$ by requiring $c_{\lambda,\lambda}=1$
(the ``monic'' choice).

The multi-variable Jacobi polynomials $\{ J_\lambda^{2/\beta,\alpha_1,
\alpha_2}(\vec{x})\}$ are orthogonal with respect to the joint
probability density \eqref{eq:2} of the \JBE\ 
\cite[Th\'eor\`eme 2]{las:pdj}.

\subsection{Multivariate hypergeometric functions}
\label{sec:mult-hyperg-funct}
Multivariate hypergeometric functions were introduced 
for general values of the Jack parameter $\sigma$ by
Kaneko \cite{kan:sia} and Kor\'anyi \cite{kor:hti}, 
generalising the definition relevant to the case $\sigma=2$
introduced by Herz \cite{her:bfo}, the Statistics applications 
of which being studied in \cite{con:snd, jam:dom, mui:aom}.
Efficient numerical implementations of multi-variable hypergeometric
functions are available \cite{koe:tee}.

They are defined as a sum over partitions as
\begin{equation}
  \label{eq:6}
  \fourIdx{}{p}{(\sigma)}{q}{F}(a_1,\ldots, a_p; b_1,\ldots, b_q;\vec{x})
\coloneq \sum_{\lambda} \frac{[a_1]_\lambda^{(\sigma)}\cdots 
[a_p]_\lambda^{(\sigma)}}{[b_1]_\lambda^{(\sigma)}\cdots 
[b_q]_\lambda^{(\sigma)}|\lambda|!} C_\lambda^{(\sigma)}(\vec{x}),
\end{equation}
where $[a]_\lambda^{(\sigma)}$ is the generalised Pochhammer symbol
defined by
\begin{equation}
  \label{eq:7}
  [a]_\lambda^{(\sigma)} \coloneq \prod_{i=1}^n \left( a-\frac{i-1}\sigma
\right)_{\lambda_i},
\end{equation}
and the classical Pochhammer symbol $(a)_n$ is
\begin{equation}
  \label{eq:148}
  (a)_n \coloneq \frac{\Gamma(a+n)}{\Gamma(a)} = \prod_{j=0}^{n-1} (a+j),
\qquad a\in\C, n\in \N.
\end{equation}
The reader familiar with hypergeometric functions of a single variable
(recapitulated in \eqref{eq:79} below) will recognise the generalisation
 \eqref{eq:6}. 

For general values of the parameters $a_1,\ldots,a_p,b_1,\ldots,b_q$
the series in \eqref{eq:6} converges absolutely for all
$\vec{x}\in\C^n$ if $p\leq q$ and for $\vec{x}$ in some ball if
$p=q+1$ \cite{kan:sia}.  However, if any of the ``upper'' parameters,
$a_1$ say, is equal to a negative integer $-m$, $m>0$, then the series
contains only finitely-many terms and defines a multi-variable
symmetric polynomial of degree $mn$.
\subsection{Some useful identities}

Yan undertook one of the first systematic studies of multi-variable 
hypergeometric functions for arbitrary $\sigma>0$ and proved a number of 
formul\ae\ and identities, including
the Pfaff-like formula \cite[eq.~(35)]{yan:aco}
\begin{equation}
\label{eq:11}
  \twoFone{\sigma}(a,b;c;\vec{x}) 
=\prod_{i=1}^n(1-x_i)^{-a} \twoFone{\sigma}
\left( a, c-b; c; \frac{-x_1}{1-x_1},\ldots, \frac{-x_n}{1-x_n}\right).
\end{equation}
(The case $\sigma=2$ was derived in \cite[Theorem 7.4.3]{mui:aom}.)

A number of integral representations are also available.  We mention
here, and will use below, the formula due to
Kaneko \cite{kan:sia}:
\begin{multline}
\label{eq:22}
  \int_0^1\cdots\int_0^1 \prod_{i=1}^n x_i^{a-1}(1-x_i)^{b-1} \prod_{j=1}^m
(x_i-t_j) |\Delta(\vec{x})|^{2/\sigma}\rmd^n\vec{x} \\
=S_n(a+m,b,1/\sigma) \fourIdx{}2{(1/\sigma)}1{F}\left( -n, 
{\sigma}({a+b+m-1}) +n-1; \sigma({a+m-1}); \vec{t}\right),
\end{multline}
valid for $\Re\{a\}>0, \Re\{b\}>0$, $\Re\{1/\sigma\}>-\min\{1/n, 
\Re\{a\}/(n-1), \Re\{b\}/(n-1)\}$, and $\vec{t}=(t_1,\ldots,t_m)$.
In \eqref{eq:22} and below we use $\rmd^{n}\vec{x}$ as a shorthand 
for $\rmd x_1\cdots \rmd x_n$.
\section{Calculations for finite-size matrices}\label{sec:main-calculations}
In this section we collect some formul\ae\ for the distribution and 
density of the smallest eigenvalue of a \JBE\ matrix of fixed finite
size $N\times N$.
\subsection{Probability distribution of the smallest eigenvalue}
If $\phi_1$ is the smallest eigenvalue of the \JBE\ then, for any
constant $\xi\in\R$,
\begin{equation}
  \label{eq:4}
 F_{\phi_1}(\xi) \coloneq \Prob(\phi_1 \leq \xi) = 1 - \Prob(\phi_1 > \xi).
\end{equation}
As all eigenvalues are between $0$ and $1$, we obviously have
\begin{equation}
  \label{eq:128}
  F_{\phi_1}(\xi) = \left\{ \begin{array}{cc}
                            0, & \xi\leq0, \\
                            1, & \xi\geq1, 
                          \end{array}\right.
\end{equation}
so it will be sufficient to find expressions for the probability
$\Prob(\phi_1 > \xi)$ in the range $0<\xi<1$.
\begin{proposition}\label{prop:finite_N}
  Let $\phi_1$ be the smallest eigenvalue of the joint distribution
\eqref{eq:2} with $\alpha_1\in\N$.  Then, for $0<\xi<1$,
\begin{equation}
  \label{eq:24}
    \Prob(\phi_1 > \xi) = (1-\xi)^{N(1+\alpha_2+(N-1)\beta/2)}
 \fourIdx{}2{(\beta/2)}1{F}\left( -N, 1-N-\frac2\beta(\alpha_2+1);
\frac{2\alpha_1}\beta;\xi \vec{1}^{\alpha_1}\right). 
\end{equation}
\end{proposition}
\dimostrazione
Recalling that $x_1,\ldots,x_n$ are un-ordered eigenvalues, we integrate
the joint probability density \eqref{eq:2}, to get
\begin{align}
  \Prob(\phi_1 > \xi) &= \Prob(x_1 > \xi, x_2>\xi, \ldots, x_N>\xi) 
\nonumber \\
&= \frac1{S_N(\alpha_1+1, \alpha_2+1, \beta/2)}  
\int_\xi^1\cdots \int_\xi^1
\prod_{i=1}^N
x_i^{\alpha_1} (1-x_i)^{\alpha_2}|\Delta(\vec{x})|^{\beta}\,\rmd^N\vec{x}.
  \label{eq:5}
\end{align}

If we make the substitution $y_i = (x_i-\xi)/(1-\xi)$, $1\leq i \leq N$,
this maps each of the integrals to an integral over $[0,1]$, and we have
\begin{multline}
\label{eq:113}
  \Prob(\phi>\xi) = \frac{(1-\xi)^{N(1+\alpha_1+\alpha_2 + (N-1)\beta/2)}}%
{S_N(\alpha_1+1,\alpha_2+1,\beta/2)} \\\times
\int_0^1\cdots\int_0^1 \prod_{i=1}^N \left( y_i + \frac{\xi}{1-\xi}
\right)^{\alpha_1}(1-y_i)^{\alpha_2} |\Delta(\vec{y})|^\beta\,\rmd^N\vec{y}.
\end{multline}
For $\alpha_1\in\N$ this integral can be evaluated by means
of Kaneko's integral \eqref{eq:22} to get
\begin{multline}
  \label{eq:23}
  \Prob(\phi>\xi) =
  (1-\xi)^{N(1+\alpha_1+\alpha_2 + (N-1)\beta/2)} \\ \times%
 \fourIdx{}2{(\beta/2)}1{F}\left( -N, \frac{2}{\beta}(\alpha_1+\alpha_2+1)+
N-1; \frac{2\alpha_1}{\beta}; \frac{-\xi}{1-\xi}\vec{1}^{\alpha_1}\right).
\end{multline}
Up to this point we have followed Borodin and Forrester's paper \cite{bor:isa}
(our \eqref{eq:23} is equation (3.16) therein).  The only novel step in
the proof is to simplify the argument of the multivariate
hypergeometric function in \eqref{eq:23} by applying the Pfaff-like
identity \eqref{eq:11} to give \eqref{eq:24}.  \finire

Based on \eqref{eq:23}, Borodin and Forrester proved the asymptotic
scaling limit \eqref{eq:150} for the smallest eigenvalue.

We have also a formula for $\Prob(\phi_1>\xi)$ in terms of 
multi-variable Jacobi polynomials.
\begin{corollary} \label{cor:Jacobi_polys}
  With $\phi_1$ and $\alpha_1\in\N$ as above, an alternative
expression for the probability in Proposition \ref{prop:finite_N}
is, for $0<\xi<1$,
\begin{equation}
  \label{eq:153}
  \Prob(\phi_1>\xi) = (1-\xi)^{N(1+\alpha_1+\alpha_2+(N-1)\beta/2)}
\frac{P_{(N^{\alpha_1})}^{\beta/2,-1+2/\beta,-1+2(\alpha_2+1)/\beta}
\left(\frac{-\xi}{1-\xi}\vec{1}^{\alpha_1}\right)}%
{P_{(N^{\alpha_1})}^{\beta/2,-1+2/\beta,-1+2(\alpha_2+1)/\beta}
(\vec{0}^{\alpha_1})},
\end{equation}
where $P_\lambda^{\sigma,a,b}(\vec{x})$ is the multi-variable Jacobi
polynomial, and an explicit expression for the denominator in
\eqref{eq:153} is
\begin{multline}
   P_{(N^{\alpha_1})}^{\beta/2,-1+2/\beta,-1+2(\alpha_2+1)/\beta}
(\vec{0}^{\alpha_1}) = \frac{(-1)^{N\alpha_1} \alpha_1! (N\alpha_1)! 
(1+2\alpha_1/\beta)_{N-1}}{(N-1)! (N\beta/2)_{\alpha_1}}
\\\times 
\prod_{i=1}^{\alpha_1} \frac1{(N-1+2(\alpha_2+1+i)/\beta)_N}.
\end{multline}
In these expressions $\vec{0}^n$ is a shorthand for $(0,\ldots,0)\in
\R^n$.
\end{corollary}
Corollary \ref{cor:Jacobi_polys} will be proved in Section
\ref{sec:main-result:2}.  We also note that explicitly computable 
recursions for coefficients in the series expansion in powers of $\xi$
for $F_{\phi_1}(\xi)$ have been derived in \cite{for:csf}.

\subsection{Probability density of the smallest eigenvalue}
Our main interest is in the probability \emph{distribution} function
of the smallest eigenvalue $\phi_1$ of the \JBE. However
with little effort we can derive a formula for a probabilty density
in terms of a multi-variable hypergeometric function,
that will also be used to prove a key differentiation identity
(Corollary \ref{cor:derivative} below).
\begin{proposition}
If $\alpha_1\in\N$, a marginal probability density function
for the smallest eigenvalue $\phi_1$ of the Jacobi
$\beta$-Ensemble \eqref{eq:2} is given 
by
  \begin{multline}
  \label{eq:13}
   Z_N(\alpha_1, \alpha_2, \beta)   \phi_1^{\alpha_1} 
(1-\phi_1)^{\alpha_2 + (N-1)(1+\alpha_2 +N\beta/2)} \\
  \times \twoFone{\beta/2} \left( 1-N, 2 - N - \frac2\beta(\alpha_2+1); 
\frac{2\alpha_1}\beta + 2 ; \phi_1\vec{1}^{\alpha_1} \right),
\end{multline}
for $0\leq\phi_1\leq1$, where the normalisation constant is 
\begin{equation}
  \label{eq:14}
  Z_N(\alpha_1, \alpha_2, \beta) \coloneq   \frac{N S_{N-1}(\alpha_1 + 1 
  + \beta, \alpha_2 +1, \beta/2)}%
  {S_N(\alpha_1 + 1, \alpha_2 + 1, \beta/2)} .
\end{equation}
\end{proposition}
\dimostrazione
The joint probability density function of the \emph{ordered} eigenvalues
of the \JBE\ is
\begin{equation}
  \label{eq:144}
\frac{N!}{S_N(\alpha_1+1, \alpha_2+1, \beta/2)}  \prod_{i=1}^N
\phi_i^{\alpha_1} (1-\phi_i)^{\alpha_2}|\Delta(\vecphi)|^{\beta},\qquad
0\leq \phi_1\leq\cdots\leq\phi_N \leq 1.
\end{equation}
This has the same functional form as \eqref{eq:2}, except for the
factor $N!$ in the numerator to account for the ordering of the
variables.  To derive the marginal density function for $\phi_1$ we
integrate out all the other variables
\begin{align}
&\frac{N!}{S_N(\alpha_1+1, \alpha_2+1, \beta/2)}
\int_{\phi_1}^1 \int_{\phi_2}^1 \cdots\int_{\phi_{N-1}}^1  \prod_{i=1}^N
\phi_i^{\alpha_1} (1-\phi_i)^{\alpha_2}|\Delta(\vecphi)|^{\beta}
  \,\rmd\phi_2\cdots\rmd\phi_N \nonumber \\
&=\frac{N!/(N-1)!}{S_N(\alpha_1+1, \alpha_2+1, \beta/2)}
\int_{\phi_1}^1 \int_{\phi_1}^1 \cdots\int_{\phi_{1}}^1  \prod_{i=1}^N
\phi_i^{\alpha_1} (1-\phi_i)^{\alpha_2}|\Delta(\vecphi)|^{\beta}
  \,\rmd\phi_2\cdots\rmd\phi_N,  \label{eq:145}
\end{align}
un-ordering the $N-1$ integrations.  With the change of variables
$y_i=(\phi_i-\phi_1)/(1-\phi_1)$ for $i=2,\ldots,N$, this multiple
integral becomes
\begin{multline}
  \label{eq:146}
  \frac{N}{S_N(\alpha_1+1, \alpha_2+1, \beta/2)} \phi_1^{\alpha_1}
(1-\phi_1)^{\alpha_2} (1-\phi_1)^{(N-1)(1+\alpha_1+\alpha_2+\beta) +
(N-1)(N-2)\beta/2} \\
\times \int_{0}^1 \int_{0}^1 \cdots\int_{0}^1  \prod_{i=2}^N
\left(y_i + \frac{\phi_1}{1-\phi_1}\right)^{\alpha_1} 
(1-y_i)^{\alpha_2}y_i^\beta |\Delta(\vec{y})|^{\beta}
  \,\rmd^{N-1}\vec{y}
\end{multline}
where, in a slightly unusual notation $\vec{y}=(y_2,\ldots,y_N)\in
\R^{N-1}$.  The $(N-1)$-fold multiple integral may be evaluated by
means of Kaneko's integral \eqref{eq:22} to give
\begin{multline}
  \label{eq:12}
  \frac{N S_{N-1}(\alpha_1 + 1 + \beta, \alpha_2 +1, \beta/2)}%
  {S_N(\alpha_1 + 1, \alpha_2 + 1, \beta/2)} 
  \phi_1^{\alpha_1} (1-\phi_1)^{\alpha_2 + (N-1)(1+\alpha_1 + \alpha_2
+N\beta/2)} \\
  \times \twoFone{\beta/2} \left( 1-N, \frac2\beta(\alpha_1+\alpha_2+1)
+ N; \frac{2\alpha_1}\beta + 2 ; \frac{-\phi_1}{1-\phi_1}\vec{1}^{\alpha_1}
 \right).
\end{multline}
By an application of the Pfaff-like identity \eqref{eq:11}, this may be
re-written as \eqref{eq:13}.  \finire

The formula \eqref{eq:13} for the probability density function was first
derived by Dumitriu \cite{dum:sed}, with a different method of proof.
Slightly different, but equivalent, multivariable hypergeometric function
representations for the probability density function have
been given in \cite{dre:cwb}.

\begin{corollary}  For $\alpha_1\in\N$   \label{cor:derivative}
  we have the derivative identity
  \begin{multline}
    \label{eq:25}
    \frac{\rmd}{\rmd \xi} \left( (1-\xi)^{N(1+\alpha_2+(N-1)\beta/2)}
 \fourIdx{}2{(\beta/2)}1{F}\left( -N, 1-N-\frac2\beta(\alpha_2+1);
\frac{2\alpha_1}\beta;\xi \vec{1}^{\alpha_1}\right)\right)
\\ = - Z_N(\alpha_1,\alpha_2,\beta) \xi^{\alpha_1}
(1-\xi)^{N(1+\alpha_2+(N-1)\beta/2)-1} \\ \times
       \twoFone{\beta/2} \left( 1-N, 2 - N - \frac2\beta(\alpha_2+1); 
\frac{2\alpha_1}\beta + 2 ; \xi\vec{1}^{\alpha_1} \right),
  \end{multline}
for all $\xi\in\C$ except possibly $\xi=1$.
\end{corollary}
\dimostrazione
From \eqref{eq:4} and \eqref{eq:24} above the probability distribution 
function of $\phi_1$, the smallest eigenvalue, is 
  \begin{equation}
F_{\phi_1}(\xi) =    1 - (1-\xi)^{N(1+\alpha_2+(N-1)\beta/2)}
 \fourIdx{}2{(\beta/2)}1{F}\left( -N, 1-N-\frac2\beta(\alpha_2+1);
\frac{2\alpha_1}\beta;\xi \vec{1}^{\alpha_1}\right),
  \end{equation}
for $0\leq\xi\leq1$,
and a probability density function is given by \eqref{eq:13}.  
The result \eqref{eq:25} follows because the density function agrees with
the derivative of the distribution function at points of continuity. By
analytic continuation the identity persists outside of the interval
$0<\xi<1$.
\finire

We remark that the result \eqref{eq:25} does not seem easy to prove in a
direct way starting from the definition \eqref{eq:6} of the
multivariate hypergeometric functions.  A similar observation was made
by Forrester \cite{for:era} who found the analogous identity at the
level of $\fourIdx{}1{(\sigma)}1{F}$ multivariate hypergeometric
functions.

Later, we will want to take the limit $N\to\infty$, so we record here
the asymptotic behaviour of $Z_N$ in this limit.
\begin{lemma}  \label{lem:Z}
  As $N\to\infty$ we have 
\begin{equation}
  \label{eq:15}
  Z_N(\alpha_1, \alpha_2, \beta) \sim \frac{ \Gamma(1+\beta/2) 
(\beta/2)^{2\alpha_1+1}}%
  {\Gamma(1+\alpha_1)\Gamma(1+\alpha_1+\beta/2)} 
N^{2(\alpha_1+1)}.
\end{equation}
\end{lemma}
\dimostrazione
  Using the value \eqref{eq:16} for the Selberg integrals, and cancelling
common factors we get
\begin{multline}
  \label{eq:17}
  Z_N(\alpha_1, \alpha_2, \beta) =  \frac{N\Gamma(1+\beta/2)}{\Gamma(1 +
N\beta/2)} \frac{\Gamma(\alpha_1 + \alpha_2 + 2 + (N-1)\beta) }%
{\Gamma(\alpha_1 + 1 + (N-1)\beta/2) \Gamma(\alpha_2+1+(N-1)\beta/2)}
\\
\times \prod_{k=0}^{N-2} \frac{\Gamma(\alpha_1 + 1 + \beta + k\beta/2)
\Gamma(\alpha_1+\alpha_2+2+(N+k-1)\beta/2)}%
 {\Gamma(\alpha_1 + 1 + k\beta/2)\Gamma(\alpha_1+\alpha_2+\beta+2+
(N+k-2)\beta/2)}.
\end{multline}
Re-writing the factors appearing in the product in \eqref{eq:17} as
\begin{equation}
  \label{eq:18}
   \frac{\Gamma(\alpha_1 + 1 + (k+2)\beta/2)
\Gamma(\alpha_1+\alpha_2+2+(N+k-1)\beta/2)}%
 {\Gamma(\alpha_1 + 1 + k\beta/2)\Gamma(\alpha_1+\alpha_2+2+
(N+k)\beta/2)}
\end{equation}
we realise that many factors cancel in the product over $k$ and we are
left with
\begin{multline}
  \label{eq:19}
  \prod_{k=0}^{N-2} \frac{\Gamma(\alpha_1 + 1 + \beta + k\beta/2)
\Gamma(\alpha_1+\alpha_2+2+(N+k-1)\beta/2)}%
 {\Gamma(\alpha_1 + 1 + k\beta/2)\Gamma(\alpha_1+\alpha_2+\beta+2+
(N+k-2)\beta/2)} \\
= \frac{\Gamma(\alpha_1+1+(N-1)\beta/2)\Gamma(\alpha_1 + 1 + N\beta/2)
\Gamma(\alpha_1 + \alpha_2 + 2 + (N-1)\beta/2)}%
  {\Gamma(\alpha+1)\Gamma( \alpha_1+1+\beta/2)
   \Gamma(\alpha_1 + \alpha_2 + 2 + (N-1)\beta)}.
\end{multline}
Reuniting the product with the prefactors in \eqref{eq:17} and further
cancellation results in
\begin{equation}
  \label{eq:20}
  Z_N(\alpha_1, \alpha_2, \beta)  = \frac{N \Gamma(1+\beta/2)\Gamma(\alpha_1
+1+N\beta/2)\Gamma(\alpha_1 + \alpha_2 + 2 + (N-1)\beta/2)}%
  {\Gamma(1+N\beta/2)\Gamma(\alpha_1+1)\Gamma(\alpha_1 + 1 + \beta/2)
\Gamma(\alpha_2+1+(N-1)\beta/2)}.
\end{equation}
The asymptotic \eqref{eq:15} follows by applying the asymptotic
formula
\begin{equation}
  \label{eq:21}
  \frac{\Gamma(a + cN)}{\Gamma(b+cN)} \sim (cN)^{a-b}
\end{equation}
to the $N$-dependent factors.\finire
\section{Two-term asymptotic formula}\label{sec:two-term-asymptotic}
Our main analytic tool is going to be a two-term asymptotic formula for the
$\fourIdx{}2{(\sigma)}1F$ multi-variable hypergeometric function,
stated below, and proved in the following subsections.
\begin{theorem}
  \label{thm:two-term}
Let $a,b,c\in\C$ and $\sigma>0$ be fixed, such that $c-(i-1)/\sigma$ is not 
a negative integer for $1\leq i\leq n$.  Then with $p_1(\vec{x})\coloneq x_1 
+ \cdots +x_n$,
\begin{multline}
  \label{eq:99}
      \fourIdx{}2{(\sigma)}1{F} \left( a-N, b-N; c; \frac1{N^2}\vec{x}\right)
= \fourIdx{}0{(\sigma)}1{F}(;c;\vec{x})  + \frac1N
(c - a - b ) \rmE_1
\left\{\fourIdx{}0{(\sigma)}1{F}(;c;\vec{x}) \right\} \\ -
\frac1Np_1(\vec{x}) \fourIdx{}0{(\sigma)}1{F}(;c;\vec{x}) 
 + \Ord\left(\frac1{N^2}\right),
\end{multline}
where the error estimate is uniform for $\vec{x}$ in compact subsets
of $\C^n$, but may depend on $a,b,c,\sigma, n$.  The operator
$\rmE_1$ in \eqref{eq:99} is defined in \eqref{eq:diff2}.
\end{theorem}

Our strategy will be to split the sum defining the 
$\fourIdx{}2{(\sigma)}1F$ multi-variable hypergeometric function as
\begin{multline}
  \label{eq:76}
  \fourIdx{}2{(\sigma)}1{F} \left( a-N, b-N; c; \frac1{N^2}\vec{x}\right)
= \sum_{|\lambda| < N^{1/3}} \frac{[a-N]_\lambda^{(\sigma)}
[b-N]_\lambda^{(\sigma)}}{N^{2|\lambda|}[c]_{\lambda}^{(\sigma)} 
|\lambda|!} C_\lambda^{(\sigma)}(\vec{x}) \\
+ \sum_{|\lambda|\geq N^{1/3}} \frac{[a-N]_\lambda^{(\sigma)}
[b-N]_\lambda^{(\sigma)}}{N^{2|\lambda|}[c]_{\lambda}^{(\sigma)} 
|\lambda|!} C_\lambda^{(\sigma)}(\vec{x}),
\end{multline}
recalling that the Jack polynomials $C_\lambda^{(\sigma)}$ are
homogeneous of order $|\lambda|$.
It will turn out that the tail terms (the second sum) do not
contribute significantly to the $N\to\infty$ limit.

\subsection{Preliminary results}
It is a consequence of Stirling's formula that
\begin{equation}
\label{eq:78}
    \frac{\Gamma(z+\alpha)}{\Gamma(z+\beta)} = z^{\alpha-\beta}\left(
1 + \frac{(\alpha-\beta)(\alpha+\beta-1)}{2z} +
\Ord\left(  \frac1{z^2}\right)\right),
\end{equation}
as $|z|\to\infty$.  We will need a form of this result with control on
how the error depends on the parameters $\alpha,\beta$. 
\begin{lemma} \label{lem:stirling}
  Suppose $\alpha\in\C$ is a quantity such that $|\alpha|^2/z=\littleo(1)$
as $z\to\infty$.  Then 
\begin{equation}
  \label{eq:26}
  \Gamma(z+\alpha) = \sqrt{2\pi} z^{z+\alpha-1/2}\rme^{-z}\left(
    1 + \frac{\alpha(\alpha-1)+1/6}{2z} + \Ord\left( \frac{|\alpha|^4+1}{z^2}
\right)\right),
\end{equation}
as $z\to\infty$ with $|\arg\{z\}|<\pi$.
\end{lemma}
\dimostrazione
By the classical Stirling formula \cite[6.1.37]{abr:hmf}
\begin{equation}
  \label{eq:27}
  \Gamma(z) =  \sqrt{2\pi} z^{z-1/2}\rme^{-z} \left( 1 + \frac1{12z}
    +\Ord\left( \frac1{z^2}\right) \right),
\end{equation}
 as $z\to\infty$ with $|\arg\{z\}|<\pi$.
Therefore,
\begin{align}
  \Gamma(z+\alpha) &= \sqrt{2\pi} (z+\alpha)^{z+\alpha-1/2}\rme^{-z-\alpha} 
\left( 1 + \frac1{12(z+\alpha)}  \nonumber
    +\Ord\left( \frac1{(z+\alpha)^2}\right) \right) \\
&= \sqrt{2\pi} (z+\alpha)^{z+\alpha-1/2}\rme^{-z-\alpha} 
\left( 1 + \frac1{12z}
    +\Ord\left( \frac{|\alpha|+1}{z^2}\right) \right),\label{eq:31}
\end{align}
using
\begin{equation}
  \label{eq:28}
  \frac1{z+\alpha} = \frac1z + \Ord\left( \frac{|\alpha|}{z^2}\right)
\end{equation}
and
\begin{equation}
  \label{eq:29}
  \frac1{(z+\alpha)^2} = \Ord\left( \frac1{z^2}\right)
\end{equation}
provided $|\alpha/z|\leq1/2$.

Now,
\begin{align}
  (z+\alpha)^{z+\alpha-1/2} &= z^{z+\alpha-1/2}\left( 1+\frac\alpha{z}
\right)^{z+\alpha-1/2} \nonumber \\
                            &= z^{z+\alpha-1/2} \exp\left( \left(
 z+\alpha-\frac12\right)\log\left( 1 + \frac\alpha{z}\right)\right) \nonumber \\
&= z^{z+\alpha-1/2} \exp\left( \left(
 z+\alpha-\frac12\right)\left( \frac\alpha{z} - \frac{\alpha^2}{2z^2}
+\Ord\left( \frac{|\alpha|^3}{z^3} \right) \right)\right) \nonumber \\
&=  z^{z+\alpha-1/2} \exp\left( \alpha + \frac{\alpha(\alpha-1)}{2z}
+\Ord \left( \frac{|\alpha|^3+|\alpha|^2}{z^2}\right) \right) \nonumber \\
&=  z^{z+\alpha-1/2} \rme^\alpha \exp\left(\frac{\alpha(\alpha-1)}{2z}\right)
\left( 1  +\Ord \left( \frac{|\alpha|^3+|\alpha|^2}{z^2}\right) \right) 
\nonumber \\
&=  z^{z+\alpha-1/2} \rme^\alpha \left(1 + \frac{\alpha(\alpha-1)}{2z}
 + \Ord \left( \frac{|\alpha|^4+|\alpha|^2}{z^2}\right) \right),
  \label{eq:30}
\end{align}
provided $|\alpha|^2/|z|\leq 1/2$.   Putting it together with \eqref{eq:31}
we get \eqref{eq:26}.   \finire
\begin{corollary}  \label{cor:gamma}
    Suppose $\alpha$ and $\beta$ satisfy $|\alpha|^2/z=\littleo(1)$
and  $|\beta|^2/z=\littleo(1)$
as $z\to\infty$.  Then 
\begin{equation}
  \label{eq:32}
  \frac{\Gamma(z+\alpha)}{\Gamma(z+\beta)} = z^{\alpha-\beta}\left(
1 + \frac{(\alpha-\beta)(\alpha+\beta-1)}{2z} +
\Ord\left(  \frac{|\alpha|^4+|\beta|^4+1}{z^2}\right)\right).
\end{equation}
\end{corollary}
\dimostrazione
Applying Lemma \ref{lem:stirling} and cancelling common factors,
\begin{align}
  \label{eq:33}
  \frac{\Gamma(z+\alpha)}{\Gamma(z+\beta)} &= z^{\alpha-\beta}\left(
1 + \frac{\alpha(\alpha-1)-1/6}{2z} + \Ord\left( \frac{|\alpha|^4 + 1}{z^2}
\right)\right) \nonumber \\ &\qquad\qquad\times\left(
1 + \frac{\beta(\beta-1)-1/6}{2z} + \Ord\left( \frac{|\beta|^4 + 1}{z^2}
\right)\right)^{-1} \nonumber \\
&=  z^{\alpha-\beta}\left(
1 + \frac{\alpha(\alpha-1)-\beta(\beta-1)}{2z} + 
\Ord\left( \frac{|\alpha|^4 +|\beta|^4 + 1}{z^2}
\right)\right),
\end{align}
by the binomial theorem.  Re-writing the second term gives \eqref{eq:32}.
\finire

We will also require some bounds on Pochhammer symbols
\eqref{eq:148}.
\begin{lemma}  \label{lem:pochhammer}
  Let $a\in\C$ be fixed, and $n,N\in\N$.
  \begin{enumerate}
  \item If $n\leq N$ then $|(a-N)_n|\leq (|a|+N-n+1)_n\leq (N+|a|)^n$;
  \item If $n>N$ then $|(a-N)_n|\leq (|a|)_n$.
  \end{enumerate}
\end{lemma}
\dimostrazione
We have that
\begin{align}
  (a-N)_n &= (a-N)(a-N+1)\cdots(a-N+n-1) \nonumber \\
  &=(-1)^n(N-a)(N-a-1)\cdots(N-a-(n-1)).\label{eq:62}
\end{align}
If $n\leq N$ then using the second line of \eqref{eq:62} and the
triangle inequality
\begin{align}
  |(a-N)_n| &= |N-a||N-1 - a|\cdots|N-(n-1) -a| \nonumber \\
 &\leq (N+|a|)(N-1+|a|)\cdots(N-(n-1)+|a|) = (|a|+N-n+1)_{n} \nonumber \\
 &\leq (N+|a|)^n.\label{eq:65}
\end{align}

If $n>N$ then re-ordering the product,
\begin{equation}
  \label{eq:63}
  (a-N)_n = a(a+1)\cdots(a+n-N-1) \times (a-1)(a-2)\cdots(a-N).
\end{equation}
However,
\begin{align}
  \label{eq:64}
\big|  (a-1)(a-2)\cdots(a-N) \big| &\leq (|a|+1)(|a|+2)\cdots(|a|+N) \nonumber
\\
& \leq (|a|+n-N)(|a|+n-N+1) \cdots (|a|+n-1),
\end{align}
so from \eqref{eq:63} we end up with
\begin{align}
  |(a-N)_n| &\leq |a|(|a|+1)\cdots(|a|+n-N-1) \nonumber \\
&\qquad \times  (|a|+n-N)(|a|+n-N+1) 
\cdots (|a|+n-1) \nonumber \\ &= (|a|)_n.
  \label{eq:66}
\end{align}
\finire

We shall test certain sums for convergence by comparison with the classical 
hypergeometric series
\begin{equation}
  \label{eq:79}
\fourIdx{}p{}q{F}(a_1,\ldots,a_p;b_1,\ldots,b_q;z) = 
  \sum_{k=0}^\infty \frac{(a_1)_k\cdots (a_p)_k}{(b_1)_k\cdots(b_q)_k}
\frac{z^k}{k!}.
\end{equation}
For generic choice of parameters $a_1,\ldots,a_p,b_1,\ldots,b_q\in\C$
the power series \eqref{eq:79} is known to have radius of convergence 
$r=1$ if $p=q+1$ and infinite radius of convergence if $p\leq q$
\cite[\S 2.2]{sla:ghf}.  The exceptions to these rules are when the series
has only a finite number of terms and reduces to a polynomial in 
$z$. This can happen when one of the parameters $a_1,\ldots,a_p$
is a negative integer.
\subsection{The main contribution}
\label{sec:main-contribution}
We start by analysing the first sum on the right-hand side of
\eqref{eq:76}.
\begin{proposition}  \label{prop:main-contribution}
  Let $a,b$ be fixed quantities.  Then
  \begin{align}
    \label{eq:83}
    &\sum_{|\lambda| < N^{1/3}} \frac{[a-N]_\lambda^{(\sigma)}
[b-N]_\lambda^{(\sigma)}}{N^{2|\lambda|}[c]_{\lambda}^{(\sigma)} 
|\lambda|!} C_\lambda^{(\sigma)}(\vec{x}) \\&= \sum_{|\lambda|<N^{1/3}}
\frac{1}{[c]_\lambda^{(\sigma)} |\lambda|!} 
\left( 1 + \frac1{N}\left(\left( \frac2\sigma
(n-1) - a - b\right)\rmE_1 - \rmD_2\right)\right) C_\lambda^{(\sigma)}(\vec{x})
+\Ord\left(\frac1{N^2}\right) \nonumber
  \end{align}
where the implied constant may depend on $n,a,b,c$ and $\sigma$ but is uniform
for $\vec{x}$ in compact subsets of $\C^n$.
\end{proposition}
Making use of the representation
\begin{equation}
  \label{eq:80}
  [a-N]_\lambda^{(\sigma)} = \prod_{i=1}^n \left( a-N-\frac{i-1}\sigma
\right)_{\lambda_i} = 
\prod_{i=1}^n (-1)^{\lambda_i} \frac{\Gamma(1+N+(i-1)/\sigma-a)}%
{\Gamma(1+N+(i-1)/\sigma-a-\lambda_i)}
\end{equation}
we can write, using Corollary \ref{cor:gamma} for the 
ratios of gamma functions,
\begin{align}
  [a-N]_\lambda^{(\sigma)} & [b-N]_\lambda^{(\sigma)} \nonumber \\ &= 
\prod_{i=1}^{n} (-1)^{2|\lambda_i|} \frac{\Gamma(1+N+(i-1)/\sigma-a)
\Gamma(1+N+(i-1)/\sigma-b)}{\Gamma(1+N+(i-1)/\sigma-a-\lambda_i)
\Gamma(1+N+(i-1)/\sigma-b-\lambda_i)} \nonumber \\
&= \prod_{i=1}^{n} N^{2\lambda_i}\left( 1 + \frac{\lambda_i}{2N}
\left( \frac{2}\sigma(i-1) +1 - 2a - \lambda_i\right) +
\Ord\left( \frac{|\lambda_i|^4+(\sigma'i)^4}{N^2}\right)\right) \nonumber \\
&\qquad\quad\times
\left( 1 + \frac{\lambda_i}{2N}
\left( \frac{2}\sigma(i-1) +1 - 2b - \lambda_i\right) +
\Ord\left( \frac{|\lambda_i|^4+(\sigma'i)^4}{N^2}\right)\right) ,
\end{align}
where $\sigma'=\max\{1,\sigma^{-1}\}$.  This leads to 
\begin{align}
  \label{eq:81}
    [&a-N]_\lambda^{(\sigma)} [b-N]_\lambda^{(\sigma)} \nonumber \\ &= 
\prod_{i=1}^{n} N^{2\lambda_i}\left( 1 + \frac{\lambda_i}N 
\left( \frac2\sigma(i-1)-a-b-\lambda_i+1\right) +
\Ord\left( \frac{|\lambda_i|^4+(\sigma'i)^4}{N^2}\right)\right)\nonumber \\
&=N^{2|\lambda|} \Bigg( 1 + \frac{|\lambda|}N\left( \frac2{\sigma}(n-1)-a-b
\right) - \frac1{N}\sum_{i=1}^n \lambda_i\left( \lambda_i - 1 - \frac{2(i-1)}
\sigma + \frac{2(n-1)}\sigma\right) \nonumber \\ &\qquad\qquad
+ \Ord\left( \frac{\|\lambda\|_4^4 + {\sigma'}^4n^5}{N^2}\right)\Bigg).
\end{align}
We adopt here a convenient shorthand $\|\lambda\|_4^4=\lambda_1^4+\cdots
\lambda_n^4$.
Recalling the actions \eqref{eq:diff3}, \eqref{eq:diff4} of the operators
$\rmE_1$, $\rmD_2$ on Jack polynomials, 
\begin{multline}
  \label{eq:82}
  [a-N]_\lambda^{(\sigma)} [b-N]_\lambda^{(\sigma)} C_{\lambda}^{(\sigma)}
(\vec{x}) =  N^{2|\lambda|} \Bigg(\left( 1 + \frac1{N}\left( \left( \frac2\sigma
(n-1) - a - b\right)\rmE_1 - \rmD_2\right)\right) C_\lambda^{(\sigma)}(\vec{x})
\\+ \Ord\left( \frac{\|\lambda\|_4^4 + {\sigma'}^4n^5}{N^2}
|C_{\lambda}^{(\sigma)}(\vec{x})|\right)\Bigg).
\end{multline}

Substituting \eqref{eq:82} for the numerator in \eqref{eq:83} leads 
\textit{via} cancellation of the factors $N^{2|\lambda|}$ to the
sum on the right-hand side of \eqref{eq:83}.
 To prove the error estimate it will be sufficient to demonstrate that
\begin{equation}
  \label{eq:67}
  \sum_{\lambda} \left|\frac{\|\lambda\|_4^4}{[c]_{\lambda}^{(\sigma)} 
|\lambda|!} C_\lambda^{(\sigma)}(\vec{x})\right| < \infty
\end{equation}
uniformly for $\vec{x}$ in compact subsets of $\C^n$.
We do this following a method elaborated by Kaneko \cite{kan:sia}. 
Namely: there
exists a constant $C$ depending only on $n$ such that
\begin{equation}
  \label{eq:68}
  |C_{\lambda}^{(\sigma)}(\vec{x})| \leq C |\lambda|^{n/2}\left( 
\sigma(\sigma'n)^{3/2} \|\vec{x}\|_{\infty} \right)^{|\lambda|},
\end{equation}
where $\|\vec{x}\|_\infty = \max\{|x_1|,\ldots,|x_n|\}$ and
$\sigma'=\max\{1,\sigma^{-1}\}$ \cite[Lemma 1]{kan:sia}.  Thus there
exists a constant $R>0$ depending only on $n, \sigma$ such that
$\|\lambda\|_4^4 |C_{\lambda}^{(\sigma)}(\vec{x})| 
\leq C R^{|\lambda|} \|\vec{x}\|_\infty^{|\lambda|}$.  We also
may observe that
\begin{equation}
  \label{eq:69}
  |\lambda|! = (\lambda_1+ \cdots+\lambda_n)! \geq \lambda_1!\cdots
\lambda_n!.
\end{equation}
So
\begin{align}
  \sum_{\lambda}\left| \frac{\|\lambda\|_4^4}{[c]_{\lambda}^{(\sigma)} 
|\lambda|!} C_\lambda^{(\sigma)}(\vec{x}) \right| &\leq
C \sum_{\lambda} \frac{R^{|\lambda|}\|\vec{x}\|_\infty^{|\lambda|}}%
 {|[c]_\lambda^{(\sigma)}||\lambda|!}  \nonumber \\
&\leq C \prod_{i=1}^n \sum_{\lambda_i=0}^\infty 
\frac{(R\|\vec{x}\|_\infty)^{\lambda_i}}{|(c-(i-1)/\sigma)_{\lambda_i}|
\lambda_i!},
  \label{eq:70}
\end{align}
using the definition \eqref{eq:7} of the generalised Pochhammer symbol
together with \eqref{eq:69}.  It can be seen that \eqref{eq:70} is
bounded by comparing each factor to a convergent $\fourIdx{}{0}{}{1}{F}$
hypergeometric series.  \finire

\subsection{Bounding the tail terms}
Using Kaneko's bound for Jack polynomials from the end of 
subsubsection \ref{sec:main-contribution} we get
\begin{equation}
   \bigg| \sum_{|\lambda| \geq N^{1/3}} 
\frac{[a-N]_\lambda^{(\sigma)}[b-N]_\lambda^{(\sigma)}
}{N^{2|\lambda|}[c]_{\lambda}^{(\sigma)} 
|\lambda|!} C_\lambda^{(\sigma)}(\vec{x}) \bigg| \leq
C \sum_{|\lambda| \geq N^{1/3}} 
\frac{|[a-N]_\lambda^{(\sigma)}[b-N]_\lambda^{(\sigma)}|
(R\|\vec{x}\|_\infty)^{|\lambda|}}{N^{2|\lambda|}|[c]_{\lambda}^{(\sigma)}| 
\lambda_1!\cdots\lambda_n!},
  \label{eq:71}
\end{equation}
for some constant $C$ depending only on $n$ and $R$ depending only on
$n$ and $\sigma$.  If $\lambda$ is a partition with $|\lambda|\geq N^{1/3}$ 
then, as $\lambda_1$ is the largest part, we must have 
$\lambda_1\geq N^{1/3}/n$.  Factorising the generalised Pochhammer
symbols according to \eqref{eq:7} we achieve the inequality
\begin{multline}
  \label{eq:72}
   \bigg| \sum_{|\lambda| \geq N^{1/3}} 
\frac{[a-N]_\lambda^{(\sigma)}[b-N]_\lambda^{(\sigma)}
}{N^{2|\lambda|}[c]_{\lambda}^{(\sigma)} 
|\lambda|!} C_\lambda^{(\sigma)}(\vec{x}) \bigg| \leq 
\\C\bigg( \sum_{\lambda_1=N^{1/3}/n}^\infty 
\frac{|(a-N)_{\lambda_1}||(b-N)_{\lambda_1}|(R\|\vec{x}\|_\infty)^{\lambda_1}}%
{N^{2\lambda_1}|(c)_{\lambda_1}| \lambda_1!}  \bigg)\\
\times \prod_{i=2}^n\sum_{\lambda_i=0}^\infty   
\frac{|(a-(i-1)/\sigma-N )_{\lambda_i}||(b-(i-1)/\sigma-N)_{\lambda_i}
|(R\|\vec{x}\|_\infty)^{\lambda_i}}%
{N^{2\lambda_i}|(c-(i-1)/\sigma)_{\lambda_i}| \lambda_i!}.
\end{multline}
\begin{proposition}  \label{prop:tail}
 For every $\rho>0$ and compact set $K\subseteq \C^n$, there exists a 
constant $C_{\rho,K}$ (which may additionally depend on $a, b, c, \sigma, n$)
 such that for every $\vec{x}\in K$, and all $N$ sufficiently large, we
have
  \begin{equation}
    \label{eq:777}
    \bigg| \sum_{|\lambda| \geq N^{1/3}} 
\frac{[a-N]_\lambda^{(\sigma)}[b-N]_\lambda^{(\sigma)}
}{N^{2|\lambda|}[c]_{\lambda}^{(\sigma)} 
|\lambda|!} C_\lambda^{(\sigma)}(\vec{x}) \bigg|\leq \frac{C_{\rho,K}}{N^\rho}
  \end{equation}
\end{proposition}
\dimostrazione
For brevity, let us define $a'=a-(i-1)/\sigma$, $b'=b-(i-1)/\sigma$,
$c'=c-(i-1)/\sigma$.
The following three estimates give us the bound we need:
\begin{enumerate}
\item  Using part 2.\ of Lemma \ref{lem:pochhammer},
  \begin{align}
    \sum_{\lambda_i>N} \frac{|(a'-N )_{\lambda_i}||(b'-N)_{\lambda_i}
|(R\|\vec{x}\|_\infty)^{\lambda_i}}%
{N^{2\lambda_i}|(c')_{\lambda_i}| \lambda_i!} &\leq
 \frac1{N^N}\sum_{\lambda_i>N} \frac{(|a'|)_{\lambda_i}(|b'|)_{\lambda_i}}%
{|(c')_{\lambda_i}|\lambda_1!} \left( \frac{R\|\vec{x}\|_\infty}{N} 
\right)^{\lambda_i}  \nonumber \\
&=\Ord\left( \frac1{N^N}\right) 
\label{eq:73}
  \end{align}
provided we additionally choose $N$ large enough so that $R\|\vec{x}\|_\infty/
N<1$, and the absolute convergence of a $\fourIdx{}2{}1{F}$ hypergeometric
series in the unit disc.
\item Using part 1.\ of Lemma \ref{lem:pochhammer},
  \begin{align}
      \sum_{\lambda_i=0}^N \frac{|(a'-N )_{\lambda_i}||(b'-N)_{\lambda_i}
|(R\|\vec{x}\|_\infty)^{\lambda_i}}%
{N^{2\lambda_i}|(c')_{\lambda_i}| \lambda_i!} &\leq
\sum_{\lambda_i=0}^N \frac{(N+|a'|)^{\lambda_i} (N+|b'|)^{\lambda_i}
(R\|\vec{x}\|_\infty)^{\lambda_i}}{N^{2\lambda_i}|(c')_{\lambda_i}| 
\lambda_i!} \nonumber \\
&=\sum_{\lambda_i=0}^N \frac{((1+|a'|/N)(1+|b'|/N)
R\|\vec{x}\|_\infty)^{\lambda_i}}{|(c')_{\lambda_i}| 
\lambda_i!} \nonumber \\
&<\infty \label{eq:74}
  \end{align}
by comparison with a $\fourIdx{}0{}1{F}$ hypergeometric series.
\item For the sum over $\lambda_1$ between $N^{1/3}/n$ and $N$ we 
again use part 1.\ of Lemma \ref{lem:pochhammer}, as in the previous step
leading to
\begin{align}
  \label{eq:75}
        \sum_{N^{1/3}/n<\lambda_1\leq N} 
&\frac{|(a-N )_{\lambda_i}||(b-N)_{\lambda_i}
|(R\|\vec{x}\|_\infty)^{\lambda_i}}%
{N^{2\lambda_i}|(c)_{\lambda_i}| \lambda_i!} \nonumber \\ &\leq
 \sum_{N^{1/3}/n<\lambda_1\leq N} \frac{((1+|a|/N)(1+|b|/N)
R\|\vec{x}\|_\infty)^{\lambda_i}}{|(c)_{\lambda_i}| 
\lambda_i!} \nonumber \\
&\leq \frac1{\Gamma(N^{1/3}/n+1)}
        \sum_{N^{1/3}/n<\lambda_1\leq N} \frac{((1+|a|/N)(1+|b|/N)
R\|\vec{x}\|_\infty)^{\lambda_i}}{|(c)_{\lambda_i}| } \nonumber \\
&=\Ord\left( \frac1{\Gamma(N^{1/3}/n+1)} \right),
\end{align}
where the series has been compared to a $\fourIdx{}1{}1{F}$ 
hypergeometric series.
\end{enumerate}
We use 1.\ and 2.\ to bound each factor for $i\geq2$ in \eqref{eq:72}
by a constant.  We then use 1.\ and 3.\ to deduce the rapid decay
in $N$ of the remaining sum over $\lambda_1$.
\finire

With essentially the same method and calculations we can bound
similarly the tail of two further series.
\begin{proposition}  \label{prop:more-tails}
 For every $\rho>0$ and compact set $K\subseteq \C^n$, there exists a 
constant $C_{\rho,K}$ (which may additionally depend on $a,b,c, \sigma, n$)
 such that for every $\vec{x}\in K$, and all $N$ sufficiently large, we
have
  \begin{equation}
    \label{eq:77}
    \bigg| \sum_{|\lambda| \geq N^{1/3}} 
\frac{1}{[c]_{\lambda}^{(\sigma)} |\lambda|!} 
\left( \left( \frac2\sigma
(n-1) - a - b\right)\rmE_1 - \rmD_2\right)
C_\lambda^{(\sigma)}(\vec{x}) \bigg|\leq \frac{C_{\rho,K}}{N^\rho}
  \end{equation}
and
  \begin{equation}
    \label{eq:77a}
    \bigg| \sum_{|\lambda| \geq N^{1/3}} 
\frac{1}{[c]_{\lambda}^{(\sigma)} |\lambda|!} 
C_\lambda^{(\sigma)}(\vec{x}) \bigg|\leq \frac{C_{\rho,K}}{N^\rho}.
  \end{equation}
\end{proposition}
\dimostrazione Using the fact that $C_{\lambda}^{(\sigma)}(\vec{x})$
is an eigenfunction of $\rmE_1$ and $\rmD_2$ with eigenvalues that
depend only polynomially on the parts of $\lambda$, we may follow the
proof of the preceding Proposition \ref{prop:tail} making only trivial
changes. \finire
\subsection{Asymptotic formula}
\begin{proposition} \label{prop:partial-result}
For fixed $a, b$, we have
\begin{multline}
  \label{eq:87}
    \fourIdx{}2{(\sigma)}1{F} \left( a-N, b-N; c; \frac1{N^2}\vec{x}\right)
= \fourIdx{}0{(\sigma)}1{F}(;c;\vec{x}) \\ + \frac1N
\left(\left( \frac2\sigma (n-1) - a - b\right) 
\rmE_1 - \rmD_2\right)\fourIdx{}0{(\sigma)}1{F}(;c;\vec{x})
 + \Ord\left(\frac1{N^2}\right),
\end{multline}
where the implied constant may depend on $n, a, b, \sigma$ but is
uniform for $\vec{x}$ in compact subsets of $\C^n$.
\end{proposition}
\dimostrazione
Starting from \eqref{eq:76} and
as a consequence of the bound \eqref{eq:777} of Proposition
\ref{prop:tail}, we have 
\begin{equation}
  \label{eq:84}
    \fourIdx{}2{(\sigma)}1{F} \left( a-N, b-N; c; \frac1{N^2}\vec{x}\right)
= \sum_{|\lambda| < N^{1/3}} \frac{[a-N]_\lambda^{(\sigma)}
[b-N]_\lambda^{(\sigma)}}{N^{2|\lambda|}[c]_{\lambda}^{(\sigma)} 
|\lambda|!} C_\lambda^{(\sigma)}(\vec{x}) 
+ \Ord\left( \frac1{N^\rho}\right)
\end{equation}
for any $\rho>0$.  By Proposition \ref{prop:main-contribution} this is
\begin{multline}
  \label{eq:85}
    \fourIdx{}2{(\sigma)}1{F} \left( a-N, b-N; c; \frac1{N^2}\vec{x}\right)
= \sum_{|\lambda|<N^{1/3}}
\frac{1}{[c]_\lambda^{(\sigma)} |\lambda|!} C_\lambda^{(\sigma)}(\vec{x})
\\+\frac1N
\sum_{|\lambda|<N^{1/3}}
\frac{1}{[c]_\lambda^{(\sigma)} |\lambda|!} 
\left(\left( \frac2\sigma
(n-1) - a - b\right)\rmE_1 - \rmD_2\right) C_\lambda^{(\sigma)}(\vec{x})
+\Ord\left(\frac1{N^2}\right).
\end{multline}
By Proposition \ref{prop:more-tails} we can complete the sums in 
\eqref{eq:85} without affecting the error estimate to get
\begin{multline}
  \label{eq:86}
    \fourIdx{}2{(\sigma)}1{F} \left( a-N, b-N; c; \frac1{N^2}\vec{x}\right)
= \sum_{\lambda}
\frac{1}{[c]_\lambda^{(\sigma)} |\lambda|!} C_\lambda^{(\sigma)}(\vec{x})
\\+\frac1N
\sum_{\lambda}
\frac{1}{[c]_\lambda^{(\sigma)} |\lambda|!} 
\left(\left( \frac2\sigma
(n-1) - a - b\right)\rmE_1 - \rmD_2\right) C_\lambda^{(\sigma)}(\vec{x})
+\Ord\left(\frac1{N^2}\right).
\end{multline}
This is \eqref{eq:87}, recognising
\begin{equation}
  \label{eq:88}
  \fourIdx{}0{(\sigma)}1{F}(;c;\vec{x}) =  \sum_{\lambda}
\frac{1}{[c]_\lambda^{(\sigma)} |\lambda|!} C_\lambda^{(\sigma)}(\vec{x})
\end{equation}
\finire

Given that $\fourIdx{}0{(\sigma)}1{F}(;c;\vec{x})$ is continuous, this
already proves 
\begin{equation}
  \label{eq:149}
\lim_{N\to\infty} \fourIdx{}2{(\sigma)}1{F} \left( a-N, b-N; c; 
\frac1{N^2}\vec{x}\right)
= \fourIdx{}0{(\sigma)}1{F}(;c;\vec{x}),
\end{equation}
locally uniformly in $\vec{x}\in\C^n$---the leading-order of Theorem
\ref{thm:two-term}.

Our final task in this section will be to put the ``$1/N$'' term of
\eqref{eq:87} into a nicer form.
\mathversion{bold}
\subsection{Partial Differential Equation satisfied by $\fourIdx{}0{(\sigma)}1{F}$}
\mathversion{normal}

It is a Theorem of Yan \cite[Theorem 2.1]{yan:aco} and Kaneko 
\cite[Theorem 2]{kan:sia} that
if $c-(i-1)/\sigma$ is not a negative integer for any $1\leq i\leq n$ then
the unique solution to system of equations
\begin{multline}
  \label{eq:89}
  x_i(1-x_i) \frac{\partial^2\Psi}{\partial x_i^2} + 
\left( c - \frac{n-1}\sigma - \left( a + b + 1 - \frac{n-1}\sigma\right)
x_i\right) \frac{\partial\Psi}{\partial x_i} \\
+\frac1{\sigma} \left( \sum_{\substack{j=1\\j\neq i}}^n \frac{x_i(1-x_i)}%
{x_i-x_j} \frac{\partial\Psi}{\partial x_i} - 
\sum_{\substack{j=1\\j\neq i}}^n \frac{x_j(1-x_j)}{x_i-x_j}
\frac{\partial\Psi}{\partial x_j}\right) = ab\Psi,
\end{multline}
$1\leq i\leq n$,
subject to $\Psi(\vec{x})$ being symmetric in its variables and analytic at
$\vec{x}=\vec{0}$, is 
\begin{equation}
  \label{eq:90}
  \Psi(\vec{x}) = \Psi(\vec{0}) \fourIdx{}2{(\sigma)}1{F}(a,b; c; \vec{x}).
\end{equation}
This result for $\sigma=2$ was first proved by Muirhead \cite{mui:sop},
having been conjectured, apparently, by A.~G.~Constantine.  Muirhead
also shows how the system \eqref{eq:89} can be degenerated to give
the holonomic system of equations for $\fourIdx{}0{(2)}{1}{F}$
multivariate hypergeometric functions, which can easily be generalised
for arbitrary $\sigma$ as follows.
\begin{proposition}
  \label{prop:0f1system}
Provided that $c-(i-1)/\sigma$ is not a negative integer for any $1\leq 
i\leq n$,
the multivariate hypergeometric function $\fourIdx{}0{(\sigma)}1{F}(c;
\vec{x})$ is the unique solution of the $n$ differential equations
\begin{equation}
  \label{eq:91}
  x_i \frac{\partial^2\Psi}{\partial x_i^2} + 
\left( c - \frac{n-1}\sigma \right) \frac{\partial\Psi}{\partial x_i} 
+\frac1{\sigma} \left( \sum_{\substack{j=1\\j\neq i}}^n \frac{x_i}%
{x_i-x_j} \frac{\partial\Psi}{\partial x_i} - 
\sum_{\substack{j=1\\j\neq i}}^n \frac{x_j}{x_i-x_j}
\frac{\partial\Psi}{\partial x_j}\right) = \Psi,
\end{equation}
$1\leq i\leq n$,
subject to the constraints that $\Psi(\vec{x})$ is symmetric in its variables, 
is analytic at $\vec{x}=\vec{0}$ and satisfies $\Psi(\vec{0})=1$.
\end{proposition}
\dimostrazione
Since we now know \eqref{eq:149} that
\begin{equation}
  \label{eq:84b}
  \fourIdx{}{0}{(\sigma)}1F (;c;\vec{x}) = 
  \lim_{N\to\infty} \fourIdx{}{2}{(\sigma)}1F\left( -N, -N; c; \frac1{N^2}
\vec{x}\right),
\end{equation}
we set $a=b=-N$ and make the change of variables $x_i\mapsto x_i/N^2$ in
\eqref{eq:89} to get
\begin{multline}
  \label{eq:92}
  N^2x_i\left(1-\frac{x_i}{N^2}\right) \frac{\partial^2\Psi}{\partial x_i^2} + 
\left( \left(c - \frac{n-1}\sigma\right)N^2 - \left( 1 - 2N - \frac{n-1}\sigma
\right) x_i\right) \frac{\partial\Psi}{\partial x_i} \\
+\frac{N^2}{\sigma} \left( \sum_{\substack{j=1\\j\neq i}}^n 
\frac{x_i(1-x_i/N^2)}%
{x_i-x_j} \frac{\partial\Psi}{\partial x_i} - 
\sum_{\substack{j=1\\j\neq i}}^n \frac{x_j(1-x_j/N^2)}{x_i-x_j}
\frac{\partial\Psi}{\partial x_j}\right) = N^2\Psi.
\end{multline}
Dividing through by $N^2$ and letting $N\to\infty$ we recover 
\eqref{eq:91}. \finire
\begin{corollary} \label{cor:pde}
Under the same condition on $c$ as in Proposition 
\ref{prop:0f1system},
  the function $\Psi(\vec{x}) = \fourIdx{}0{(\sigma)}1{F}(;c;\vec{x})$ is
a solution to the partial differential equation
\begin{equation}
  \label{eq:93}
  \rmD_2 \Psi = p_1(\vec{x})\Psi - \left( c- \frac{2}{\sigma}(n-1)\right)
\rmE_1\Psi,
\end{equation}
where $p_1(\vec{x})=x_1+\cdots+x_n$.
\end{corollary}
\dimostrazione
We multiply through the $i$th equation \eqref{eq:91}, satisfied by
$\Psi(\vec{x}) = \fourIdx{}0{(\sigma)}1{F}(;c;\vec{x})$, by $x_i$:
\begin{equation}
  \label{eq:94}
  x_i^2 \frac{\partial^2\Psi}{\partial x_i^2} + 
\left( c - \frac{n-1}\sigma \right)x_i \frac{\partial\Psi}{\partial x_i} 
+\frac1{\sigma} \left( \sum_{\substack{j=1\\j\neq i}}^n \frac{x_i^2}%
{x_i-x_j} \frac{\partial\Psi}{\partial x_i} - 
\sum_{\substack{j=1\\j\neq i}}^n \frac{x_ix_j}{x_i-x_j}
\frac{\partial\Psi}{\partial x_j}\right) = x_i\Psi.
\end{equation}
Since
\begin{equation}
  \label{eq:95}
  \frac{x_ix_j}{x_i-x_j} = x_j + \frac{x_j^2}{x_i-x_j},
\end{equation}
\eqref{eq:94} is equivalent to
\begin{multline}
  \label{eq:96}
  x_i^2 \frac{\partial^2\Psi}{\partial x_i^2} + 
\left( c - \frac{n-1}\sigma \right)x_i \frac{\partial\Psi}{\partial x_i} 
\\+\frac1{\sigma} \left( \sum_{\substack{j=1\\j\neq i}}^n \frac{x_i^2}%
{x_i-x_j} \frac{\partial\Psi}{\partial x_i} - 
\sum_{\substack{j=1\\j\neq i}}^n x_j\frac{\partial\Psi}{\partial x_j} -
\sum_{\substack{j=1\\j\neq i}}^n \frac{x_j^2}{x_i-x_j}
\frac{\partial\Psi}{\partial x_j}\right) = x_i\Psi.
\end{multline}
Summing over \eqref{eq:96} for $i=1,\ldots,n$, we arrive at
\begin{equation}
  \label{eq:97}
  \rmD_2 \Psi + \left( c- \frac{n-1}\sigma\right) \rmE_1 \Psi 
- \frac{n-1}{\sigma} \rmE_1 \Psi = p_1(\vec{x}) \Psi.
\end{equation}
Re-arranged, this is \eqref{eq:93}. \finire

\dimostrazionea{Theorem \ref{thm:two-term}}
We use Corollary \ref{cor:pde} to replace the term $\rmD_2\left\{
\fourIdx{}{0}{(\sigma)}1{F}(;c;\vec{x})\right\}$ in \eqref{eq:87} by
\begin{equation}
  \label{eq:98}
  p_1(\vec{x}) \fourIdx{}{0}{(\sigma)}1{F}(;c;\vec{x}) +
\left(\frac{2}{\sigma}(n-1)-c\right) \rmE_1 \left\{
\fourIdx{}{0}{(\sigma)}1{F}(;c;\vec{x})\right\}.
\end{equation}
The resultant cancellation of terms involving $\sigma$ leads to
\eqref{eq:99}.  \finire

\section{Main Results}\label{sec:main-result}
\subsection{Proof of Theorem \ref{thm:main}}
Looking-back to \eqref{eq:24} we had
\begin{multline}
  \label{eq:101}
      \Prob(N^2\phi_1 > x) = \left(1-\frac{x}{N^2}\right)^{N(1+\alpha_2+
(N-1)\beta/2)}\\\times
 \fourIdx{}2{(\beta/2)}1{F}\left( -N, 1-N-\frac2\beta(\alpha_2+1);
\frac{2\alpha_1}\beta;\frac{x}{N^2} \vec{1}^{\alpha_1}\right). 
\end{multline}
Standard asymptotic arguments give
\begin{equation}
  \label{eq:102}
  \left(1-\frac{x}{N^2}\right)^{N(1+\alpha_2+
(N-1)\beta/2)} = \rme^{-\beta x/2} \left( 1 -
\left( 1 + \alpha_2 - \frac{\beta}2\right)\frac{x}N +
\Ord\left( \frac1{N^2} \right)\right)
\end{equation}
uniformly for ${x}$ in compact sets.

In the result of Theorem \ref{thm:two-term},
taking $\vec{x}$ to be a constant multiple of $\vec{1}^n$ we have
\begin{multline}
  \label{eq:100}
      \fourIdx{}2{(\sigma)}1{F} \left( a-N, b-N; c; \frac{x}{N^2}\vec{1}^n
\right)
= \fourIdx{}0{(\sigma)}1{F}(;c;x\vec{1}^n)  + \frac1N
(c - a - b ) x\frac{\rmd}{\rmd x}
\left\{\fourIdx{}0{(\sigma)}1{F}(;c;x\vec{1}^n) \right\} \\ -
\frac{n}N x\, \fourIdx{}0{(\sigma)}1{F}(;c;x\vec{1}^n) 
 + \Ord\left(\frac1{N^2}\right),
\end{multline}
which may be applied in \eqref{eq:101} with $a=0$, $b=1-2(\alpha_2+1)/\beta$, 
$c=2\alpha_1/\beta$, $\sigma=\beta/2$ and $n=\alpha_1$
(so that $c-(i-1)/\sigma = 2(\alpha_1+1-i)/\beta$ is not a negative
integer for $i=1,\ldots,\alpha_1$, justifying the use of Proposition
\ref{prop:0f1system}), to give
\begin{multline}
  \label{eq:104}
   \fourIdx{}2{(\beta/2)}1{F}\left( -N, 1-N-\frac2\beta(\alpha_2+1);
\frac{2\alpha_1}\beta;\frac{x}{N^2} \vec{1}^{\alpha_1}\right)
\\= \fourIdx{}0{(\beta/2)}1F\left(
;\frac{2\alpha_1}\beta; x\vec{1}^{\alpha_1} \right)  
+\frac1N \left( \frac{2\alpha_1}\beta + \frac2\beta(\alpha_2 +1) -1\right)
x\frac{\rmd}{\rmd x}\left\{ \fourIdx{}0{(\beta/2)}1F\left(
;\frac{2\alpha_1}\beta; x\vec{1}^{\alpha_1} \right) \right\} 
\\- \frac{\alpha_1}N x\, \fourIdx{}0{(\beta/2)}1F\left(
;\frac{2\alpha_1}\beta; x\vec{1}^{\alpha_1} \right)  + \Ord\left(\frac1{N^2}
\right).
\end{multline}

Combining with \eqref{eq:102} we get
\begin{align}
    \Prob(N^2\phi_1 > x) &= \rme^{-\beta x/2} \fourIdx{}0{(\beta/2)}1F\left(
;\frac{2\alpha_1}\beta; x\vec{1}^{\alpha_1} \right) \nonumber \\
&\qquad + \frac1N\left( \frac2\beta(\alpha_1 + \alpha_2 + 1) - 1 \right) x
\rme^{-\beta x/2} \frac{\rmd}{\rmd x}\left\{ \fourIdx{}0{(\beta/2)}1F\left(
;\frac{2\alpha_1}\beta; x\vec{1}^{\alpha_1} \right) \right\} 
\nonumber \\ & \qquad-\frac1N\left( 1 + \alpha_1 + \alpha_2 - \frac\beta2
\right) x\rme^{-\beta x/2} \fourIdx{}0{(\beta/2)}1F\left(
;\frac{2\alpha_1}\beta; x\vec{1}^{\alpha_1} \right)  + \Ord\left( \frac1{N^2}
\right) \nonumber \\
 &= \rme^{-\beta x/2} \fourIdx{}0{(\beta/2)}1F\left(
;\frac{2\alpha_1}\beta; x\vec{1}^{\alpha_1} \right) \nonumber \\
&\qquad + \frac{x}N\left(\frac{2}\beta(\alpha_1 + \alpha_2 +1)-1\right) 
\frac{\rmd}{\rmd x} \left\{ \rme^{-\beta x/2}\fourIdx{}0{(\beta/2)}1F\left(
;\frac{2\alpha_1}\beta; x\vec{1}^{\alpha_1} \right)  \right\}  \nonumber \\
&\qquad + \Ord\left( \frac1{N^2} \right).
  \label{eq:103}
\end{align}

Setting $\xi = x/N^2$ in \eqref{eq:25} of Corollary \ref{cor:derivative}
  \begin{multline}
\label{eq:105}
 N^2   \frac{\rmd}{\rmd x} \left( \left(1-\frac{x}{N^2}\right)^{N(1+\alpha_2
+(N-1)\beta/2)}
 \fourIdx{}2{(\beta/2)}1{F}\left( -N, 1-N-\frac2\beta(\alpha_2+1);
\frac{2\alpha_1}\beta;\frac{x}{N^2} \vec{1}^{\alpha_1}\right)\right)
\\ = - \frac{Z_N(\alpha_1,\alpha_2,\beta)}{N^{2\alpha_1}} x^{\alpha_1}
\left(1-\frac{x}{N^2}\right)^{N(1+\alpha_2+(N-1)\beta/2)-1} \\ \times
       \twoFone{\beta/2} \left( 1-N, 2 - N - \frac2\beta(\alpha_2+1); 
\frac{2\alpha_1}\beta + 2 ; \frac{x}{N^2}\vec{1}^{\alpha_1} \right),
\end{multline}
where $Z_N$ was defined in \eqref{eq:14}.  A further application of
(the leading order of) Theorem \ref{thm:two-term} and equation \eqref{eq:102},
and Lemma \ref{lem:Z} for the asymptotics of $Z_N$, brings \eqref{eq:105}
to
\begin{multline}
\label{eq:106}
  \frac{\rmd}{\rmd x}\left( \rme^{-\beta x/2} \fourIdx{}0{(\beta/2)}1F\left(
;\frac{2\alpha_1}\beta; x\vec{1}^{\alpha_1} \right) \right) \\=
\frac{-\Gamma(1+\beta/2)}{\Gamma(1+\alpha_1)\Gamma(1+\alpha_1+\beta/2)}
\left( \frac\beta2 \right)^{2\alpha_1+1} x^{\alpha_1} 
\rme^{-\beta x/2}
\fourIdx{}0{(\beta/2)}1F\left(
;\frac{2\alpha_1}\beta+2; x\vec{1}^{\alpha_1} \right).
\end{multline}
We use \eqref{eq:106} to remove the derivative term from \eqref{eq:103},
giving
\begin{multline}
  \label{eq:107}
    \Prob(N^2\phi_1 > x) = \rme^{-\beta x/2} \fourIdx{}0{(\beta/2)}1F\left(
;\frac{2\alpha_1}\beta; x\vec{1}^{\alpha_1} \right) \\
 - \frac{x^{1+\alpha_1}}N\left((\alpha_1 + \alpha_2 +1)-\frac\beta2\right) 
\left(\frac\beta2\right)^{2\alpha_1} 
\frac{\Gamma(1+\beta/2)}{\Gamma(1+\alpha_1)\Gamma(1+\alpha_1+\beta/2)}
\\\times \rme^{-\beta x/2}
\fourIdx{}0{(\beta/2)}1F\left(
;\frac{2\alpha_1}\beta+2; x\vec{1}^{\alpha_1} \right) +
\Ord\left(\frac1{N^2}\right).
\end{multline}
This completes the proof. \finire

That the first-order correction term in \eqref{eq:107} is proportional
to the derivative of the leading term implies that the finite-$N$ behaviour
can be interpreted, up to an error of order $\Ord(N^{-2})$ as a correction
to the width: if we let
\begin{equation}
  \label{eq:158}
  F_\infty(x)\coloneq 1 -  \rme^{-\beta x/2} \fourIdx{}0{(\beta/2)}1F\left(
;\frac{2\alpha_1}\beta; x\vec{1}^{\alpha_1} \right)
\end{equation}
then we may interpret Theorem \ref{thm:main} as saying
\begin{equation}
  \label{eq:159}
  F_{N^2\phi_1}(x) = F_\infty(x) + \frac{x}N\left(\frac2\beta(\alpha_1 +
\alpha_2 + 1)-1\right) F_\infty'(x) + \Ord\left(\frac1{N^2}\right).
\end{equation}
By Taylor's theorem, this is equivalent to the re-centring
\begin{equation}
  \label{eq:160}
  F_{N^2\phi_1}(x) = F_\infty\left( x\left( 1 + \frac1N\left(
\frac2\beta(\alpha_1 +
\alpha_2 + 1)-1\right) \right)\right) + \Ord\left(\frac1{N^2}\right).
\end{equation}
\subsection{Connection with Jacobi polynomials}  \label{sec:main-result:2}
We now prove the formula \eqref{eq:153} from 
Corollary \ref{cor:Jacobi_polys} giving a formula for the distribution
of the smallest \JBE\ eigenvalue in terms of multi-variable Jacobi
polynomials.

\dimostrazionea{Corollary \ref{cor:Jacobi_polys}}
We have for $\vec{x}\in\C^n$
\begin{multline}
  \label{eq:39}
  \fourIdx{}2{(\sigma)}1{F}(-N, N+1+a+b+(n-1)/\sigma; a + 1 +(n-1)/\sigma;
\vec{x}) \\= 
\frac{[-N]_{(N^n)}^{(\sigma)}[N+1+a+b+(n-1)/\sigma]_{(N^n)}^{(\sigma)}}{[a+1
+(n-1)/\sigma]_{(N^n)}^{(\sigma)} (nN)!}
P_{(N^n)}^{\sigma,a,b}(\vec{x}).
\end{multline}
This is essentially Th\'eor\`eme 5 of \cite{las:pdj}, incorporating our
different choice of normalisation of the Jacobi polynomials.  It may be proved
by observing that both sides of \eqref{eq:39} are multivariate
symmetric polynomials that satisfy the same 
partial differential equation (see \cite{yan:aco} or \cite{kan:sia} for
the PDE satisfied by $\fourIdx{}2{(\sigma)}1{F}$) and that the
leading term on both sides is proportional to $C_{(N^n)}^{(\sigma)}(\vec{x})$
with identical constant.

Comparing to \eqref{eq:23} we find the appropriate parameters for the
multivariate Jacobi polynomial are
\begin{equation}
  \label{eq:40}
  \sigma = \frac{\beta}2,\qquad a= -1 + \frac2\beta,\qquad 
b = -1 + \frac{2(\alpha_2+1)}\beta.
\end{equation}
From \eqref{eq:23} and the identity \eqref{eq:39} with parameters as above
we have
\begin{align}
  \label{eq:41}
  \Prob(&\phi_1 > \xi) = (1-\xi)^{N(1+\alpha_1+\alpha_2+(N-1)\beta/2)}\\
&\times
\frac{[-N]_{(N^{\alpha_1})}^{(\beta/2)}[N-1+2(\alpha_1+\alpha_2+1)/\beta
]_{(N^{\alpha_1})}^{(\beta/2)}}{[2\alpha_1/\beta]_{(N^{\alpha_1})}^{(\beta/2)}
 (N\alpha_1)!} P_{(N^{\alpha_1})}^{\beta/2,-1+2/\beta,-1+2(\alpha_2+1)/\beta}
\left(\frac{-\xi}{1-\xi}\vec{1}^{\alpha_1}\right).\nonumber
\end{align}
Taking the limit $\xi\to0^+$ we need to have $\Prob(\phi_1>0)=1$ and so
\begin{equation}
  \label{eq:42}
   P_{(N^{\alpha_1})}^{\beta/2,-1+2/\beta,-1+2(\alpha_2+1)/\beta}
(\vec{0}^{\alpha_1}) = \frac{[2\alpha_1/\beta]_{(N^{\alpha_1})}^{(\beta/2)}
 (N\alpha_1)!}{[-N]_{(N^{\alpha_1})}^{(\beta/2)}[N-1+2(\alpha_1+\alpha_2+
1)/\beta]_{(N^{\alpha_1})}^{(\beta/2)}}.
\end{equation}

Some of the quantities in \eqref{eq:42} simplify:  we have, starting
from \eqref{eq:7},
\begin{align}
  [-N]_{(N^n)}^{(\sigma)} &= \prod_{i=1}^n\prod_{j=1}^N \left( -N + j -1
-\frac{i-1}\sigma \right) \nonumber \\
&=\prod_{i=1}^n\prod_{j=1}^N \left( -j -\frac{i-1}\sigma \right) \nonumber \\
&=(-1)^{Nn}\prod_{i=0}^{n-1}\prod_{j=1}^{N} \left( j + \frac{i}\sigma \right),
\label{eq:46}
\end{align}
and
\begin{align}
  \left[ \frac{n}\sigma+\theta \right]_{(N^n)}^{(\sigma)} &= 
\prod_{i=1}^n\prod_{j=1}^N \left( \theta + j - 1  + \frac{n-i+1}\sigma 
\right) \nonumber \\
&=\prod_{i=1}^n\prod_{j=0}^{N-1} \left( \theta + j + \frac{i}\sigma 
\right),
  \label{eq:47}
\end{align}
so that
\begin{align}
  \frac{[n/\sigma]_{(N^n)}^{(\sigma)}}{[-N]_{(N^n)}^{(\sigma)}} &=
(-1)^{nN} \frac{\prod_{i=1}^n \prod_{j=0}^{N-1}(j+i/\sigma)}%
{\prod_{i=0}^{n-1}\prod_{j=1}^N (j+i/\sigma)} \nonumber \\
&= (-1)^{Nn} \frac{\prod_{i=1}^n (i/\sigma) \prod_{j=1}^{N-1} (j+n/\sigma)}%
{\prod_{i=0}^{n-1} (N+i/\sigma)\prod_{j=1}^{N-1}j} \nonumber \\
&= \frac{(-1)^{nN}\, n!}{(N-1)!\sigma^n} \frac{\prod_{j=1}^{N-1} (j+n/\sigma)}%
{\prod_{i=0}^{n-1} \sigma^{-1}(N\sigma+i)} \nonumber \\
&= \frac{(-1)^{nN}\, n! (1+n/\sigma)_{N-1}}{(N-1)!(N\sigma)_n}. 
  \label{eq:48}
\end{align}
In \eqref{eq:42} the ratio $[2\alpha_1/\beta]_{(N^{\alpha_1})}^{(\beta/2)}/
[-N]_{(N^{\alpha_1})}^{(\beta/2)}$ is of the form \eqref{eq:48} and
the factor $[N-1+2(\alpha_1+\alpha_2+1)/\beta]_{(N^n)}^{(\beta/2)}$ is
of the form \eqref{eq:47}.  Thus we get
\begin{multline}
  \label{eq:49}
   P_{(N^{\alpha_1})}^{\beta/2,-1+2/\beta,-1+2(\alpha_2+1)/\beta}
(\vec{0}^{\alpha_1}) = \frac{(-1)^{N\alpha_1} \alpha_1! (N\alpha_1)! 
(1+2\alpha_1/\beta)_{N-1}}{(N-1)! (N\beta/2)_{\alpha_1}}
\\\times 
\prod_{i=1}^{\alpha_1} \frac1{(N-1+2(\alpha_2+1+i)/\beta)_N}.
\end{multline}
\finire
\subsection{Numerical simulations}
In order to visualise our results better, we present in this subsection
the results of some numerical simulations, and compare them against our
theoretical predictions.

In Figures \ref{fig:eins} and \ref{fig:zwei} we present numerical results
for two different values of $\beta$: $\beta=1$ (Figure \ref{fig:eins}) and
$\beta=3$ (Figure \ref{fig:zwei}).  In both cases we plot the empirical 
distribution function for 1\,000 samples of the scaled smallest 
eigenvalue of \JBE\ random matrices for $N=30$ and $N=20$ respectively,
together with the same calculation at $N=1000$.  We observe that for
the smaller value of $N$, the empirical distribution fits well to the
two-term asymptotic prediction \eqref{eq:140} of Theorem \ref{thm:main}.
At $N=1000$, the data follows the leading order term of \eqref{eq:140}
(which is the formula \eqref{eq:150} derived by Borodin and Forrester 
\cite{bor:isa}), the order $1/N$ corrections being negligble at such 
matrix size.
\begin{figure}
  \centering
      \begin{tikzpicture}
      \begin{axis}[height=7cm,
        xmin=-0.5,
        xmax=10.5,
        xlabel={Scaled smallest eigenvalue},
        legend style={anchor=south east, at={(0.98,0.1)}},
        legend entries={{{$N\to\infty$} prediction},
                         {$N=30$ correction},
                         {$N=30$} numerics,
                        {{$N=10^3$} numerics}}]
        \addplot[dashed, thick, domain=0:10, samples=100] 
                   {1 - exp(-0.5*x)};
       \addplot[thick, domain=0:10, samples=100] 
        {1 - exp(-0.5*x) + (3.5/30 * x * exp(-0.5*x))};
       \addplot[only marks, mark=o] table {data3.dat};
      \end{axis}
      \node at (-0.75,5.25) {a)};
      \begin{scope}[xshift = 8cm]
      \begin{axis}[height=7cm,
        xmin=-0.5,
        xmax=10.5,
        xlabel={Scaled smallest eigenvalue},
        legend style={anchor=south east, at={(0.98,0.1)}},
        legend entries={{{$N\to\infty$} prediction},
                         {$N=30$ correction},
                        {{$N=10^3$} numerics}}]
        \addplot[dashed, thick, domain=0:10, samples=100] 
                   {1 - exp(-0.5*x)};
       \addplot[thick, domain=0:10, samples=100] 
        {1 - exp(-0.5*x) + (3.5/30 * x * exp(-0.5*x))};
       \addplot[only marks, mark=square] table {data4.dat};
      \end{axis}
      \node at (-0.75,5.25) {b)};
      \end{scope}
    \end{tikzpicture}
    \caption{ A plot of the empirical cumulative probability
      distribution (points) for the scaled smallest eigenvalue
      $N^2\phi_1$ for 1000 samples from the Jacobi Orthogonal Ensemble
      (i.e.\ $\beta=1$) with parameters $\alpha_1=0$, $\alpha_2=3$,
      for (a) $N=30$; (b) $N=1000$.  Also plotted are the leading term
      and the $N=30$ first-order correction term of the prediction
      \eqref{eq:123} (lines).  }
  \label{fig:eins}
\end{figure}
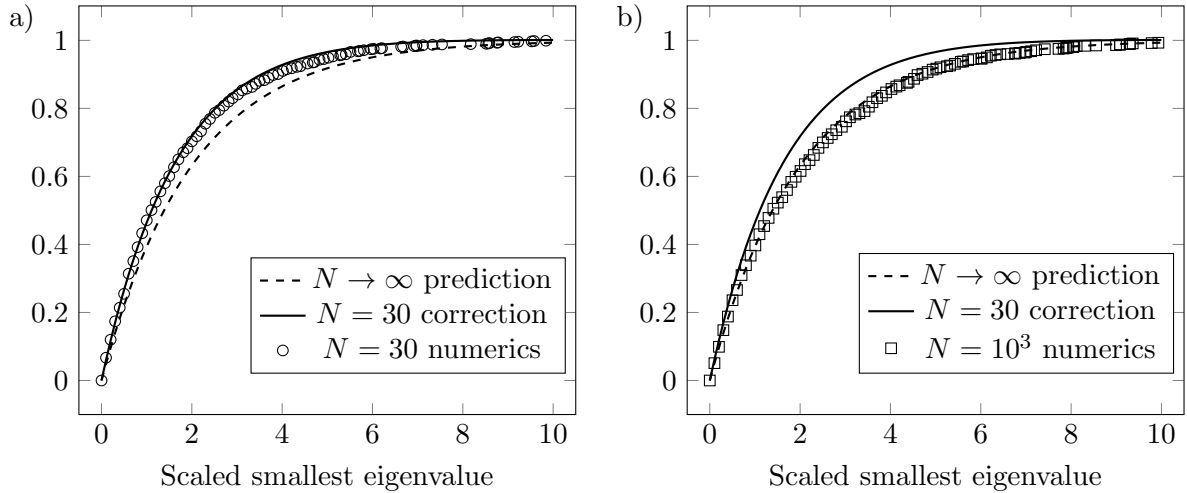
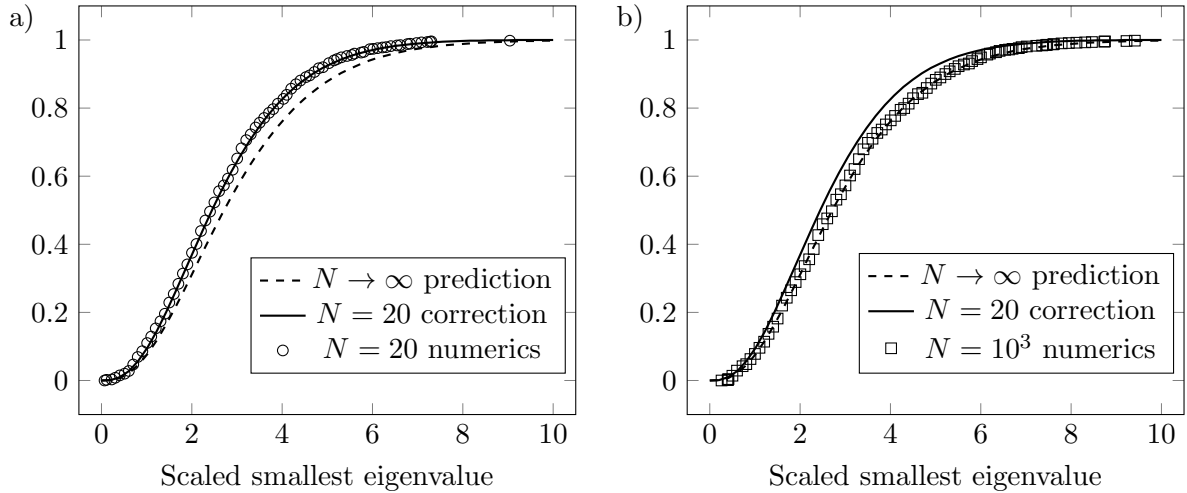
\begin{figure}
  \centering
    \begin{tikzpicture}
      \begin{axis}[height=7cm,
        xmin=-0.5,
        xmax=10.5,
        xlabel={Scaled smallest eigenvalue},
        legend style={anchor=south east, at={(0.98,0.1)}},
        legend entries={{{$N\to\infty$} prediction},
                         {$N=20$ correction},
                         {$N=20$} numerics,
                        {{$N=10^3$} numerics}}]
        \addplot[dashed, thick] table [x index=0, y index=1] {data5.dat};
        \addplot[thick] table [x index=0, y index=2] {data5.dat};
       \addplot[only marks, mark=o] table {data6.dat};
      \end{axis}
      \node at (-0.75,5.25) {a)};
      \begin{scope}[xshift = 8cm]
      \begin{axis}[height=7cm,
        xmin=-0.5,
        xmax=10.5,
        xlabel={Scaled smallest eigenvalue},
        legend style={anchor=south east, at={(0.98,0.1)}},
        legend entries={{{$N\to\infty$} prediction},
                         {$N=20$ correction},
                        {{$N=10^3$} numerics}}]
        \addplot[dashed, thick] table [x index=0, y index=1] {data5.dat};
        \addplot[thick] table [x index=0, y index=2] {data5.dat};
       \addplot[only marks, mark=square] table {data7.dat};
      \end{axis}
      \node at (-0.75,5.25) {b)};
      \end{scope}
  \end{tikzpicture}
  \caption{ A plot of the empirical cumulative probability
    distribution (points) for the scaled smallest eigenvalue
    $N^2\phi_1$ for 1000 samples from the Jacobi $\beta$ Ensemble with
    $\beta=3$, with parameters $\alpha_1=2$, $\alpha_2=1.7$, for (a)
    $N=20$; (b) $N=1000$.  Also plotted are the leading term and the
    $N=20$ first-order correction term of the prediction
    \eqref{eq:140} (lines).}
  \label{fig:zwei}
\end{figure}

The implementation details differ somewhat between Figures \ref{fig:eins}
and \ref{fig:zwei}.  To construct random matrices from the Jacobi
Orthogonal Ensemble ($\beta=1$) to generate the data for Figure
\ref{fig:eins}, we used the double-Wishart matrix construction starting
from independent $(N+1)\times N$ and $(N+7)\times N$ matrices as
described in the Introduction, which results in samples from the 
$\alpha_1=0, \alpha_2=3$ JOE.  For the theoretical prediction with
$\alpha_1=0$ we were able to use the simpler formula \eqref{eq:123}
from Section \ref{sec:small_alpha} below, rather than \eqref{eq:140}.

To generate the $\beta=3$ data for Figure \ref{fig:zwei} there is no
double-Wishart construction available, so we instead implemented the
algorithm of Killip and Nenciu \cite{kil:mmf} which can generate
samples from the \JBE\ for arbitrary $\beta>0$.  For calculating the
values of the multivariate hypergeometric functions needed for the
theoretical prediction \eqref{eq:140} we used the numerical routines
of Koev and Edelman \cite{koe:tee}.
  
\section{Explicit formul\ae}\label{sec:explicit-formulas}
In certain situations we are able to derive expressions for
$F_{\phi_1}$ and asymptotics for $F_{N^2\phi_1}$ that are more
explicit, and these are expounded in the present Section.

To begin with we focus on the case $\beta=2$ corresponding to the
Jacobi Unitary Ensemble.  In this case we benefit from the fact that
the multi-variable Jacobi polynomials enjoy a determinantal 
structure.  In Section \ref{sec:small_alpha} we record some formul\ae\
(for arbitrary $\beta>0$) for the special cases $\alpha_1=0$ and 
$\alpha_1=1$.

\subsection{Determinantal identities}
In order to state the determinantal identities
let us denote by $p_m^{a,b}(x)$ the $m$th monic Jacobi polynomial
orthogonal with respect to the measure $x^{a}(1-x)^{b}$
on the interval $[0,1]$.  In terms of the definition $P_m^{(a,b)}$
of Jacobi polynomials given by Szeg\H{o} \cite[Ch.~IV]{sze:op}, our 
definition satisfies
\begin{equation}
  \label{eq:34}
  p_m^{a,b}(x) = \frac{m!}{(m+a+b+1)_m} P_m^{(b,a)}(%
2x-1).
\end{equation}
As a hypergeometric function,
\begin{equation}
  \label{eq:35}
  p_m^{a,b}(x) = \frac{(-1)^mm!}{(m+a+b+1)_m} 
\binom{m+a}m \fourIdx{}2{}1{F}(-m, m+a+b+1; a+1; x).
\end{equation}

We also need the hook-length $h_\lambda$ for a partition, defined by
\begin{equation}
  \label{eq:154}
  h_\lambda \coloneq \frac{\prod_{i=1}^{\ell(\lambda)} (\lambda_i+
\ell(\lambda)-i)!}{\prod_{i<j}\lambda_i - \lambda_j - i + j},
\end{equation}
with $\ell(\lambda)$ the number of non-zero parts.

If $\sigma=1$ the following Lemma gives alternative expresions for
the multivariate 
Jacobi polynomials.
\begin{lemma}\label{lem:2f1}
  Let $\vec{x}\in\C^n$.  Then
\begin{equation}
  \label{eq:36}
  P_\lambda^{1,a,b}(\vec{x}) =
\frac{(-1)^{\lfloor n/2\rfloor}|\lambda|!}{h_\lambda}
\frac{\det(p_{\lambda_i+n-i}^{a,b}(x_j))_{i,j=1}^n}{\Delta(\vec{x})},
\end{equation}
where $h_\lambda$ is the 
hook-length of the partition $\lambda$ and $\Delta(\vec{x})$ is
the Vandermonde determinant \eqref{eq:1}.  If, 
furthermore $\vec{x}=x\vec{1}^n$, $x\in\C$, then we have
\begin{equation}
  \label{eq:52}
    P_\lambda^{1,a,b}(x\vec{1}^n) =
\frac{(-1)^{\lfloor n/2\rfloor}|\lambda|!}{h_\lambda \prod_{j=1}^{n-1}j!}
\det\left( \frac{(\lambda_i+n-i)!}{(\lambda_i+n-i-j+1)!} 
p_{\lambda_i+n-i-j+1}^{a+j-1,b+j-1}(x)  \right)_{i,j=1}^n.
\end{equation}
\end{lemma}
\dimostrazione
That the multivariate Jacobi polynomial at $\sigma=1$ has a determinant
evaluation in terms of univariate 
Jacobi polynomials is known since \cite[Th\'eor\`eme 10]{las:pdj}:
\begin{equation}\label{eq:43}
  P_\lambda^{1,a,b}(\vec{x}) = \text{const}
\frac{\det(p_{\lambda_i+n-i}^{a,b}(x_j))_{i,j=1}^n}{\Delta(\vec{x})}.
\end{equation}
Lasalle uses a different version of the Jacobi polynomials to us which 
changes the numerical value of the constant in \eqref{eq:43}.  
We can fix the constant
by observing that since our Jacobi polynomials are monic,
\begin{align}
  \frac{\det(p_{\lambda_i+n-i}^{a,b}(x_j))}{\Delta(\vec{x})}
&=\frac{\det(x^{\lambda_i+n-i}(x_j))}{\Delta(\vec{x})} +
\text{lower order terms} \nonumber \\
&= (-1)^{\lfloor n/2\rfloor}
\mathfrak{s}_\lambda(\vec{x}) + \text{lower order terms},
  \label{eq:37}
\end{align}
where $\mathfrak{s}_\lambda$ is a Schur polynomial.  The Jack polynomials
at $\sigma=1$ are proportional to Schur polynomials \cite{jac:aco}, and 
in fact
\begin{equation}
  \label{eq:38}
  C_\lambda^{(1)}(\vec{x}) = \frac{|\lambda|!}%
{h_\lambda} \mathfrak{s}_\lambda(\vec{x}).
\end{equation}
 The implication is
\begin{equation}
  \label{eq:36bis}
  P_\lambda^{1,a,b}(\vec{x}) =
\frac{(-1)^{\lfloor n/2\rfloor}|\lambda|!}{h_\lambda}
\frac{\det(p_{\lambda_i+n-i}^{a,b}(x_j))_{i,j=1}^n}{\Delta(\vec{x})}.
\end{equation}

In our applications $\vec{x}$ is a scalar multiple of $\vec{1}^n$.
To take the confluent limit $\vec{x}\to x\vec{1}^n$, we use the formula
\begin{equation}
  \label{eq:44}
  \lim_{\vec{x}\to x\vec{1}^n} \frac{\det(\varphi_i(x_j))}{\Delta(\vec{x})}
= \frac1{\prod_{j=1}^{n-1} j!} \mathcal{W}(\varphi_1,\ldots,\varphi_n)(x)
\end{equation}
proved in \cite[Lemma A.1]{ari:lew}, where $\mathcal{W}$ denotes the
Wronskian 
\begin{equation}
  \label{eq:45}
  \mathcal{W}(\varphi_1,\ldots,\varphi_n)(x)
= \det\left( \frac{\rmd^{j-1}}{\rmd x^{j-1}}
\varphi_i(x)\right)_{i,j=1}^n,
\end{equation}
and the functions $\varphi_1,\ldots,\varphi_n$ must be regular at
$x$. (A version of \eqref{eq:44} valid for polynomials
$\varphi_1,\ldots,\varphi_n$ was proved in \cite[Theorem 1]{gra:eop},
which would suffice to handle \eqref{eq:36}, but later we will have
reason to apply \eqref{eq:44} to non-polynomial functions.)
In the application we presently have in mind, $\varphi_i$ is the Jacobi 
polynomial $p_{\lambda_i+n-i}^{a,b}$ and using the fact that
\begin{equation}
  \label{eq:50}
  \frac{\rmd}{\rmd x} p_n^{a,b}(x) = n p_{n-1}^{a+1,b+1}(x),
\end{equation}
which extends to
\begin{align}
    \frac{\rmd^{j-1}}{\rmd x^{j-1}} p_n^{a,b}(x) &= n(n-1)\cdots
(n-j+2) p_{n-j+1}^{a+j-1,b+j-1}(x) \nonumber \\
&= \frac{n!}{(n-j+1)!} p_{n-j+1}^{a+j-1,b+j-1}(x),
  \label{eq:51}
\end{align}
we have 
\begin{equation}
  \label{eq:52bis}
  \mathcal{W}(p_{\lambda_1+n-1}^{a,b},\ldots,p_{\lambda_n}^{a,b})(x)
= \det\left( \frac{(\lambda_i+n-i)!}{(\lambda_i+n-i-j+1)!} 
p_{\lambda_i+n-i-j+1}^{a+j-1,b+j-1}(x)\right)_{i,j=1}^n.
\end{equation}
Combining \eqref{eq:36}, \eqref{eq:44} and \eqref{eq:52bis} we get
\eqref{eq:52}. \finire
\begin{corollary}
  For $\vec{x}\in\C^n$, and fixed $c$ such that $c-i$ is not a negative
 integer for $i=0,\ldots,n-1$, we have
\begin{equation}
 \label{eq:125} 
  \fourIdx{}0{(1)}1{F}(;c;\vec{x})
=(-1)^{\lfloor n/2\rfloor} \frac{\prod_{i=1}^n \Gamma(c+1-i) x_i^{n-c/2}}
{\Delta(\vec{x})}
\det\left( x_j^{-i/2} I_{c+i-2n}(2\sqrt{x_j})\right)_{i,j=1}^n,
\end{equation}
where $I_\nu(z)$ is the $I$-Bessel function.
If $\vec{x}=x\vec{1}^n$, $x\in\C$ we have a further
formula:
\begin{equation}
  \fourIdx{}0{(1)}1F(;c;x\vec{1}^n)
= x^{n(n-c)/2} \frac{\prod_{i=1}^n \Gamma(c+1-i)}{\prod_{j=1}^{n-1} j!}
\det\left( I_{c-n+j-i}(2\sqrt{x})\right)_{i,j=1}^n.
  \label{eq:138}
\end{equation}
\end{corollary}
A special case (with $c=2n$) of  \eqref{eq:138}
was given in \cite{for:bca}, where it was proved by using the
relationships with Painlev\'e functions \cite{for:dpe}.  We further
remark that \eqref{eq:125} could be proved in a less direct way
through the use of tau functions of hypergeometric type
\cite{orl:nsm}.

\medskip
\dimostrazione
We know that if $c-(i-1)/\sigma$ is not a negative integer for any $i=1,
\ldots,n$ then
\begin{equation}
  \label{eq:60}
  \fourIdx{}2{(\sigma)}1{F}\left(-N, -N; c; \frac1{N^2}\vec{x}\right) \to
\fourIdx{}0{(\sigma)}1{F}(;c;\vec{x})\qquad\text{as $N\to\infty$,}
\end{equation}
uniformly for $\vec{x}$ in compact subsets of $\mathbb{C}^n$.  
From \eqref{eq:39}
\begin{equation}
  \fourIdx{}2{(\sigma)}1{F}\left(-N, -N; c; \frac1{N^2}\vec{x}\right) =
  \frac{\left( [-N]^{(\sigma)}_{(N^n)}\right)^2}{[c]_{(N^n)}^{(\sigma)} (Nn)!} 
P_{(N^n)}^{\sigma,c-1-(n-1)/\sigma,-c-2N}\left( \frac1{N^2}\vec{x}\right).
  \label{eq:61}
\end{equation}
Setting $\sigma=1$ and applying \eqref{eq:36},
\begin{align}
  \fourIdx{}2{(1)}1{F}\left(-N, -N; c; \frac1{N^2}\vec{x}\right) &=
  \frac{\left( [-N]^{(1)}_{(N^n)}\right)^2}{[c]_{(N^n)}^{(1)} (Nn)!} 
P_{(N^n)}^{1,c-n,-c-2N}\left( \frac1{N^2}\vec{x}\right)
  \label{eq:108} \\
&= (-1)^{\lfloor n/2\rfloor} 
  \frac{\left( [-N]^{(1)}_{(N^n)}\right)^2}{[c]_{(N^n)}^{(1)}}
\frac{N^{n(n-1)}}{h_{(N^n)}} \frac{\det(p_{N+n-i}^{c-n,-c-2N}
(x_j/N^2))_{i,j=1}^n}{\Delta(\vec{x})}.
\nonumber 
\end{align}
It is problematic to take the limit $N\to\infty$ here directly.  The
determinant of Jacobi polynomials tends to $0$ as $N\to\infty$ but to
find the rate and leading-term some further manipulations are necessary.
These involve repeated use of the contiguous identity
\begin{equation}
  \label{eq:109}
  p_m^{a,b}(x) + \frac{m(m+b)}{(2m+a+b-1)(2m+a+b)} p_{m-1}^{a,b}(x) =
p_m^{a-1,b}(x),
\end{equation}
to add
successively to each row a multiple of the row below, in a 
recursive fashion, to get that
\begin{equation}
  \label{eq:110}
  \det\big(p_{N+n-i}^{c-n,-c-2N}
(x_j/N^2)\big)_{i,j=1}^n = \det\big(p_{N+n-i}^{c-2n+i,-c-2N}
(x_j/N^2)\big)_{i,j=1}^n,
\end{equation}
and we have
\begin{equation}
  \label{eq:111}
  \fourIdx{}2{(1)}1{F}\left(-N, -N; c; \frac1{N^2}\vec{x}\right) =
 (-1)^{\lfloor n/2\rfloor}
  \frac{\left( [-N]^{(1)}_{(N^n)}\right)^2}{[c]_{(N^n)}^{(1)}}
\frac{N^{n(n-1)}}{h_{(N^n)}} \frac{\det(p_{N+n-i}^{c-2n+i,-c-2N}
(x_j/N^2))}{\Delta(\vec{x})}.
\end{equation}
We return to the hypergeometric function representation of Jacobi polynomials
\eqref{eq:35}, 
\begin{multline}
  \label{eq:112}
  p_{N+n-i}^{c-2n+i,-c-2N}(x_j/N^2) = 
\frac{(-1)^{N+n-i}(N+n-i)!}{(1-n-N)_{N+n-i}} 
\binom{N+c-n}{N+n-i} \\
\times \fourIdx{}2{}1{F}(i-n-N, 1-n-N ;1+c-2n+i; x_j/N^2)
\end{multline}
and observe that
\begin{equation}
  \label{eq:115}
  (1-n-N)_{N+n-i} = (-1)^{N+n-i}\frac{(N+n-1)!}{(i-1)!}
\end{equation}
so
\begin{multline}
  \label{eq:117}
  p_{N+n-i}^{c-2n+i,-c-2N}(x_j/N^2) =\\
 \frac{(i-1)!\Gamma(N+c+1-n)}{\Gamma(N+n)
\Gamma(1+c-2n+i)} \fourIdx{}2{}1{F}(i-n-N, 1-n-N ;1+c-2n+i; x_j/N^2).
\end{multline}
Using \eqref{eq:78} and the fact that 
\begin{equation}
  \label{eq:118}
  \lim_{N\to\infty} 
\fourIdx{}2{}1{F}(i-n-N, 1-n-N ;1+c-2n+i; x_j/N^2) =
\fourIdx{}0{}1{F}(;1+c-2n+i; x_j) 
\end{equation}
we have the asymptotic behaviour
\begin{equation}
  p_{N+n-i}^{c-2n+i,-c-2N}(x_j/N^2) 
\sim \frac{(i-1)!N^{c+1-2n}}{\Gamma(1+c-2n+i)}
\fourIdx{}0{}1{F}(;1+c-2n+i; x_j) 
\end{equation}
as $N\to\infty$.  Putting this into the determinant from \eqref{eq:111},
\begin{equation}
  \label{eq:121}
  \det(p_{N+n-i}^{c-2n+i,-c-2N} (x_j/N^2))
\sim N^{cn+n-2n^2}\prod_{i=1}^n(i-1)! \det\left(
\frac{\fourIdx{}0{}1{F}(;1+c-2n+i; x_j)}{\Gamma(1+c-2n+i)}\right),
\end{equation}
as $N\to\infty$.

We already know (equation \eqref{eq:46}) that
\begin{equation}
  \label{eq:119}
  [-N]_{(N^n)}^{(1)} = (-1)^{nN} \prod_{i=1}^n \frac{\Gamma(N+
i)}{\Gamma(i)},
\end{equation}
and we have
\begin{equation}
  \label{eq:54}
  h_{(N^n)} = \prod_{i=1}^n\prod_{j=1}^N (i+j-1) 
  = \prod_{i=1}^n \frac{\Gamma(N+i)}{\Gamma(i)},
\end{equation}
so
\begin{align}
   \frac{\left( [-N]^{(1)}_{(N^n)}\right)^2}{[c]_{(N^n)}^{(1)} h_{(N^n)}}
&= \prod_{i=1}^n \frac{\Gamma(N+i)\Gamma(c+1-i)}{\Gamma(i)\Gamma(c+1-i+N)} 
\nonumber \\
&\sim \prod_{i=1}^n \frac{\Gamma(c+1-i)}{\Gamma(i)}N^{-c-1+2i} \nonumber \\
&= N^{-cn+n^2} \prod_{i=1}^n \frac{\Gamma(c+1-i)}{\Gamma(i)}\quad
\text{as $N\to\infty$.}
  \label{eq:120}
\end{align}
Putting \eqref{eq:119} and \eqref{eq:120} together into the right-hand
side of \eqref{eq:108} we find that all the factors of $N$ cancel,
and we recover the $N\to\infty$ limit, which yields
\begin{align}
\fourIdx{}0{(1)}1{F}(;c;\vec{x}) &= \lim_{N\to\infty}
    \fourIdx{}2{(1)}1{F}\left(-N, -N; c; \frac1{N^2}\vec{x}\right) 
\nonumber \\
&= (-1)^{\lfloor n/2\rfloor} \frac{\prod_{i=1}^n \Gamma(c+1-i)}{\Delta(\vec{x})}
\det\left(
\frac{\fourIdx{}0{}1{F}(;1+c-2n+i; x_j)}{\Gamma(1+c-2n+i)}\right)_{i,j=1}^n
  \label{eq:122}
\end{align}
We can derive an expression involving the more familar Bessel functions
by means of the identity \cite[\S7.8, eq.~(1)]{erd:htfII}
\begin{equation}
  \label{eq:124}
  \fourIdx{}0{}1{F}(;\alpha+1;x) = \Gamma(\alpha+1)x^{-\alpha/2}I_\alpha(2
\sqrt{x}),
\end{equation}
Inserting this into
\eqref{eq:122}, we have
\begin{equation}
  \fourIdx{}0{(1)}1{F}(;c;\vec{x}) = 
(-1)^{\lfloor n/2\rfloor} \frac{\prod_{i=1}^n \Gamma(c+1-i)}{\Delta(\vec{x})}
\det\left( x_j^{n-(c+i)/2} I_{c+i-2n}(2\sqrt{x_j})\right),
\end{equation}
which is equivalent to \eqref{eq:125}.

To take the confluent limit $\vec{x}\to x\vec{1}^n$, we prefer to work
with \eqref{eq:122}.  Using, for a second time, identity \eqref{eq:44},
\begin{multline}
  \label{eq:133}
  \lim_{\vec{x}\to x\vec{1}^n} \fourIdx{}0{(1)}1F(;c ;\vec{x}) =
(-1)^{\lfloor n/2 \rfloor} \frac{\prod_{i=1}^n \Gamma(c+1-i)}{\prod_{j=1}^{n-1}
  j!}\\\times
\det\left( \frac1{\Gamma(1+c-2n+i)} \frac{\rmd^{j-1}}{\rmd x^{j-1}}
\left( \fourIdx{}0{}1F (;1+c-2n+i; x)\right)\right).
\end{multline}
Since
\begin{equation}
  \label{eq:134}
  \frac{\rmd}{\rmd x} \fourIdx{}0{}1{F}(;c;x) = \frac1c 
\fourIdx{}0{}1{F}(;c+1;x),
\end{equation}
we have
\begin{equation}
  \label{eq:135}
  \frac{\rmd^{j-1}}{\rmd x^{j-1}} \fourIdx{}0{}1{F}(;c;x) = 
\frac{\Gamma(c)}{\Gamma(c+j-1)}  
\fourIdx{}0{}1{F}(;c+j-1;x),
\end{equation}
and
\begin{equation}
  \label{eq:136}
  \frac{\rmd^{j-1}}{\rmd x^{j-1}} \left( \fourIdx{}0{}1{F}(;1+c-2n+i;x)\right)
= \frac{\Gamma(1+c-2n+i)}{\Gamma(c-2n+i+j)} 
\fourIdx{}0{}1{F}(;c-2n+i+j;x).
\end{equation}
Putting this into \eqref{eq:133},
\begin{align}
\fourIdx{}0{(1)}1F(;c;x\vec{1}^n) &= (-1)^{\lfloor n/2 \rfloor} 
 \frac{\prod_{i=1}^n \Gamma(c+1-i)}{\prod_{j=1}^{n-1} j!}
\det\left( \frac{\fourIdx{}0{}1F (;c-2n+i+j; x)}%
{\Gamma(c-2n+i+j)} \right)_{i,j=1}^n \nonumber \\
&= \frac{\prod_{i=1}^n \Gamma(c+1-i)}{\prod_{j=1}^{n-1} j!}
\det\left( \frac{\fourIdx{}0{}1F (;c-n+j-i; x)}%
{\Gamma(c-n+j-i)} \right)_{i,j=1}^n,
  \label{eq:137}
\end{align}
reversing the order of the rows in the determinant.  Using \eqref{eq:124} this
becomes
\begin{equation}
  \fourIdx{}0{(1)}1F(;c;x\vec{1}^n)
 = \frac{\prod_{i=1}^n \Gamma(c+1-i)}{\prod_{j=1}^{n-1} j!}  
\det\left( x^{(-c+n-j+i)/2} I_{c-n+j-i}(2\sqrt{x})\right)_{i,j=1}^n
\end{equation}
As a final step we use the identity $\det(\theta^{i-j}a_{ij}) =
\det(a_{ij})$ to reduce this to \eqref{eq:138}.
\finire

\subsection{Smallest eigenvalue of the Jacobi Unitary Ensemble}
\begin{proposition}
  Let $F_{\phi_1}$ be the probability distribution function of
the smallest eigenvalue of the $N\times N$ Jacobi $\beta$-Ensemble
with $\beta=2$, $\alpha_1\in\N_0$ and $\alpha_2>-1$.  For
$0<\xi<1$,
\begin{multline}
  \label{eq:157}
  F_{\phi_1}(\xi) = 1 - (-1)^{N\alpha_1}
\prod_{i=1}^{\alpha_1}(N+\alpha_2+i)_N
(1-\xi)^{N(\alpha_1+\alpha_2+N)}\\
            \times    \det \left( \frac1{(N+i-j)!}p_{N+i-j}^{j-1,\alpha_2+j-1}
               \left( \frac{-\xi}{1-\xi}\right) \right)_{i,j=1}^{\alpha_1}.
\end{multline}
As $N\to\infty$ we have, for $x>0$,
\begin{equation}
  \label{eq:143bis}
F_{N^2\phi_1}(x) = 1 - \rme^{-x} \det(I_{j-i}(2\sqrt{x})) 
+ \frac{x}N
(\alpha_1 + \alpha_2) \rme^{-x}
\det(I_{2+j-i}(2\sqrt{x})) 
 + \Ord\left(\frac1{N^2}\right),
\end{equation}
where the determinants in \eqref{eq:143bis} are of size $\alpha_1\times
\alpha_1$.
\end{proposition}
\dimostrazione
In view of \eqref{eq:153} we need to specialise \eqref{eq:52} to the
rectangular partition $\lambda=(N^n)$, to get
\begin{equation}
  \label{eq:53}
    P_{(N^n)}^{1,a,b}(x\vec{1}^n) =
\frac{(-1)^{\lfloor n/2\rfloor}(nN)!}{h_{(N^n)} \prod_{j=1}^{n-1}j!}
\det\left( \frac{(N+n-i)!}{(N+n-i-j+1)!} 
p_{N+n-i-j+1}^{a+j-1,b+j-1}(x)  \right)_{i,j=1}^n.
\end{equation}
We have already seen (equation \eqref{eq:54}) that
\begin{equation}
  \label{eq:156}
h_{(N^n)} = \prod_{i=1}^n \frac{\Gamma(N+i)}{\Gamma(i)},  
\end{equation}
and, after cancellation, \eqref{eq:53} becomes
\begin{equation}
  \label{eq:55}
    P_{(N^n)}^{1,a,b}(x\vec{1}^n) =
\frac{(-1)^{\lfloor n/2\rfloor}(nN)!}{\prod_{i=1}^{n}(N+i-1)!}
\det\left( \frac{(N+n-i)!}{(N+n-i-j+1)!} 
p_{N+n-i-j+1}^{a+j-1,b+j-1}(x)  \right)_{i,j=1}^n.
\end{equation}
In the determinant in \eqref{eq:55} there is a factor $(N+n-i)!$ multiplying
the $i$th row.  If we extract these factors, they cancel the factorials in
the denominator.  Finally, we reverse the order of the rows producing
a factor $(-1)^{\lfloor n/2 \rfloor}$ and the end result
\begin{equation}
  \label{eq:56}
    P_{(N^n)}^{1,a,b}(x\vec{1}^n) = (nN)!
\det\left( \frac{1}{(N+i-j)!} 
p_{N+i-j}^{a+j-1,b+j-1}(x)  \right)_{i,j=1}^n.
\end{equation}

So, if $\beta=2$,
\begin{align}
  \label{eq:57}
  \Prob(\phi_1&>\xi) = \frac{(1-\xi)^{N(\alpha_1+\alpha_2+N)}}%
 {P_{(N^{\alpha_1})}^{1,0,\alpha_2}(\vec{0}^{\alpha_1})} 
  P_{(N^{\alpha_1})}^{1,0,\alpha_2}\left(\frac{-\xi}{1-\xi}\vec{1}^{\alpha_1}
 \right) \\
&=\frac{(N\alpha_1)!}{P_{(N^{\alpha_1})}^{1,0,\alpha_2}(\vec{0}^{\alpha_1})} 
(1-\xi)^{N(\alpha_1+\alpha_2+N)}
                        \det \left( \frac1{(N+i-j)!}p_{N+i-j}^{j-1,\alpha_2+j-1}
               \left( \frac{-\xi}{1-\xi}\right) \right)_{i,j=1}^{\alpha_1},
\nonumber
\end{align}
where, from \eqref{eq:49},
\begin{equation}
  \label{eq:58}
     P_{(N^{\alpha_1})}^{1,0,\alpha_2}
(\vec{0}^{\alpha_1}) = \frac{(-1)^{N\alpha_1} \alpha_1! (N\alpha_1)! 
(1+\alpha_1)_{N-1}}{(N-1)! (N)_{\alpha_1}}
\prod_{i=1}^{\alpha_1} \frac1{(N+\alpha_2+i))_N}.
\end{equation}
Note that since $\alpha_1!(1+\alpha_1)_{N-1} = (N+\alpha_1-1)!$ and
$(N-1)!(N)_{\alpha_1}=(N+\alpha_1-1)!$ too we get a simplified formula
\begin{equation}
  \label{eq:59}
     P_{(N^{\alpha_1})}^{1,0,\alpha_2}
(\vec{0}^{\alpha_1}) = {(-1)^{N\alpha_1} (N\alpha_1)! }
\prod_{i=1}^{\alpha_1} \frac1{(N+\alpha_2+i)_N}.
\end{equation}
Equations \eqref{eq:57} and \eqref{eq:59} yield \eqref{eq:157}.

We now turn to the two-term asymptotic formula \eqref{eq:143bis}.
With $\beta=2$ in \eqref{eq:140},
\begin{align}
\label{eq:141}
F_{N^2\phi_1}(x) &= 1 - \rme^{-x} \fourIdx{}0{(1)}1F\left(
;\alpha_1; x\vec{1}^{\alpha_1} \right) \\ \nonumber
&\qquad + \frac{x^{1+\alpha_1}}N(\alpha_1 + \alpha_2) 
\frac{\rme^{-x}}{\Gamma(1+\alpha_1)\Gamma(2+\alpha_1)}
\fourIdx{}0{(1)}1F\left(
;\alpha_1+2; x\vec{1}^{\alpha_1} \right) +
\Ord\left(\frac1{N^2}\right).
\end{align}
Using the representation \eqref{eq:138} for the multi-variable hypergeometric
functions this becomes
\begin{align}
\label{eq:142}
F_{N^2\phi_1}(x) &= 1 - \rme^{-x} \prod_{i=1}^{\alpha_1}\frac{\Gamma(\alpha_1
+1 - i)}{\Gamma(i)} \det(I_{j-i}(2\sqrt{x})) 
 \\ \nonumber
&\qquad + \frac{x}N
\frac{(\alpha_1 + \alpha_2) \rme^{-x}}{\Gamma(1+\alpha_1)\Gamma(2+\alpha_1)}
\prod_{i=1}^{\alpha_1} \frac{\Gamma(\alpha_1+3-i)}{\Gamma(i)}
\det(I_{2+j-i}(2\sqrt{x})) 
 + \Ord\left(\frac1{N^2}\right),
\end{align}
and, upon cancellation of the gamma function factors,
\begin{equation}
  \label{eq:143}
F_{N^2\phi_1}(x) = 1 - \rme^{-x} \det(I_{j-i}(2\sqrt{x})) 
+ \frac{x}N
(\alpha_1 + \alpha_2) \rme^{-x}
\det(I_{2+j-i}(2\sqrt{x})) 
 + \Ord\left(\frac1{N^2}\right).
\end{equation}
\finire

The formula \eqref{eq:143} proves a conjecture made in \cite{mor:eed}.
(The result is also implicit in the recent work \cite{for:roc}.)
It seems likely that explicit formul\ae\ along the lines of 
\eqref{eq:143} will also be available in the other privileged cases
$\beta=1$ and $\beta=4$.  The details will appear elsewhere.

\mathversion{bold}
\subsection{Small values of $\alpha_1$}\label{sec:small_alpha}
\mathversion{normal}
If $\alpha_1=0$, then most of the analysis above is quite unnecessary
and we already recover from \eqref{eq:113}
\begin{align}
\nonumber
  \Prob(\phi_1>\xi) &= \frac{(1-\xi)^{N(1+\alpha_2 + (N-1)\beta/2)}}%
{S_N(1,\alpha_2+1,\beta/2)}
\int_0^1\cdots\int_0^1 \prod_{i=1}^N 
(1-y_i)^{\alpha_2} |\Delta(\vec{y})|^\beta\,\rmd^N\vec{y} \\
&=(1-\xi)^{N(1+\alpha_2 + (N-1)\beta/2)}, \qquad 0<\xi<1,
\label{eq:114}
\end{align}
recognising the value of the Selberg integral (or
observing that both sides must be unity as $\xi\to0^+$).  So, for
$\alpha_1=0$ and $x>0$,
\begin{align}
\nonumber
  \Prob(N^2\phi_1 \leq x) &= 1 - \Prob\left(\phi_1 > \frac{x}{N^2}\right) \\
\nonumber
&= 1 - \left( 1 - \frac{x}{N^2} \right)^{N(1+\alpha_2 + (N-1)\beta/2)} \\
&= 1 - \rme^{-\beta x/2} + 
\left( 1 + \alpha_2 - \frac{\beta}2\right)\frac{x\rme^{-\beta x/2}}N +
\Ord\left( \frac1{N^2} \right),
  \label{eq:123}
\end{align}
referring-back to \eqref{eq:102}.  This generalises, to arbitrary 
$\beta>0$, Corollary 1 of \cite{mor:eed}.

If $\alpha_1=1$ then the multivariate hypergeometric funtions of 
$\alpha_1$ arguments become ordinary one-variable hypergeometric
functions.  In this case, for $0\leq\xi\leq1$,
\begin{equation}
  \label{eq:126}
  \Prob(\phi_1\leq \xi) = 1- (1-\xi)^{N(1+\alpha_2+(N-1)\beta/2)}
\fourIdx{}2{}{1}F\left( -N, 1-N-\frac2\beta(\alpha_2+1);
\frac2\beta; \xi\right),
\end{equation}
by \eqref{eq:24}.  Using Corollary~\ref{cor:Jacobi_polys}, and the fact
that multi-variable Jacobi polynomials of a single variable coincide
with classical single-variable ones, this may be expressed further
as
\begin{multline}
  \label{eq:161}  \Prob(\phi_1\leq \xi)
= 1 - (-1)^N(1-\xi)^{N(2+\alpha_2+(N-1)\beta/2)}\\\times
\frac{(N+2(\alpha_2+2)/\beta-1)_N}{(2/\beta)_N} p_N^{2/\beta-1,
2(\alpha_2+1)/\beta-1}\left( \frac{-\xi}{1-\xi} \right), \qquad 0< \xi<1.
\end{multline}

We apply Theorem \ref{thm:two-term} with $n=1$ to \eqref{eq:126} to find that
\begin{align}
  \Prob(N^2\phi_1\leq x) &= 1- \left(1-\frac{x}{N^2}\right)^{N(1+\alpha_2+
(N-1)\beta/2)}
\fourIdx{}2{}{1}F\left( -N, 1-N-\frac2\beta(\alpha_2+1);
\frac2\beta; \frac{x}{N^2} \right) \nonumber  \\
&=1 - \rme^{-\beta x/2} \left( 1 - \left( 1 + \alpha_2 -
\frac{\beta}2\right)\frac xN + \Ord\left( \frac1{N^2} \right)\right) 
\nonumber \\
&\qquad \times \bigg( \fourIdx{}0{}1F \left(;\frac2\beta; x\right) +
\frac1N\left(\frac2\beta(\alpha_2+2)-1\right) x\frac{\rmd}{\rmd x}\left(
\fourIdx{}0{}1F \left(;\frac2\beta; x\right)\right) \nonumber \\
&\qquad\qquad - \frac1N x\, \fourIdx{}0{}1F \left(;\frac2\beta; x\right) +
\Ord\left(\frac1{N^2}\right) \bigg) \nonumber \\
&= 1 - \rme^{-\beta x/2} \fourIdx{}0{}1F \left(;\frac2\beta; x\right)
+\frac{x}N\rme^{-\beta x/2} \left( 2 + \alpha_2 - \frac\beta2\right)
\bigg(\fourIdx{}0{}1F \left(;\frac2\beta; x\right) \nonumber \\
&\qquad\qquad  - 
\fourIdx{}0{}1F \left(;\frac2\beta+1; x\right)\bigg) +
\Ord\left( \frac1{N^2}\right),
  \label{eq:127}
\end{align}
using \eqref{eq:134} for the derivative of the hypergeometric function.
We may further use \eqref{eq:124} to replace hypergeometric functions
with Bessel functions, yielding
\begin{align}
  \Prob(N^2\phi_1\leq x) &= 1 - \rme^{-\beta x/2} \Gamma\left(\frac2\beta
\right)x^{1/2-1/\beta} I_{2/\beta-1}(2\sqrt{x})  \nonumber \\
&\qquad+\frac{x}N\rme^{-\beta x/2} \left( 2 + \alpha_2 - \frac\beta2\right)
\bigg(\Gamma\left( \frac2{\beta}\right) x^{1/2-1/\beta}I_{2/\beta-1}(
2\sqrt{x}) \nonumber \\
&\qquad\qquad  - \Gamma\left( \frac2\beta+1\right) x^{-1/\beta} I_{2/\beta}(
2\sqrt{x}) \bigg)
+ \Ord\left( \frac1{N^2}\right) \nonumber \\
&= 1 - \rme^{-\beta x/2} \Gamma\left(\frac2\beta
\right)x^{1/2-1/\beta} I_{2/\beta-1}(2\sqrt{x})  \nonumber \\
&\qquad+\frac{x}N\rme^{-\beta x/2} \left( 2 + \alpha_2 - \frac\beta2\right)
\Gamma\left( \frac2{\beta}\right) x^{-1/\beta} \bigg(x^{1/2}I_{2/\beta-1}(
2\sqrt{x}) \nonumber \\
&\qquad\qquad  - \frac2\beta I_{2/\beta}(
2\sqrt{x}) \bigg)
+ \Ord\left( \frac1{N^2}\right).
  \label{eq:129}
\end{align}
We note that 
\begin{equation}
  \label{eq:130}
  \sqrt{x} I_{2/\beta-1}(2\sqrt{x}) = \frac2\beta I_{2/\beta}(2\sqrt{x})
+ \sqrt{x} I_{2/\beta+1}(2\sqrt{x}):
\end{equation}
an application of the Bessel function identity \cite[\S7.11,
eq.~(23)]{erd:htfII}
\begin{equation}
  \label{eq:131}
  I_\nu(z) = \frac{z}{2\nu}\left( I_{\nu-1}(z) - I_{\nu+1}(z)\right).
\end{equation}
Putting \eqref{eq:130} into \eqref{eq:129} gives
\begin{multline}
  \label{eq:132}
  \Prob(N^2\phi_1\leq x) = 1 - \rme^{-\beta x/2} \Gamma\left(\frac2\beta
\right)x^{1/2-1/\beta} I_{2/\beta-1}(2\sqrt{x}) 
 \\ +\frac{x^{3/2-1/\beta}}N \rme^{-\beta x/2} \left( 2 + \alpha_2 - 
\frac\beta2\right)
\Gamma\left( \frac2{\beta}\right) I_{2/\beta+1}(2\sqrt{x})
+ \Ord\left( \frac1{N^2}\right).
\end{multline}
This result is consistent with \eqref{eq:143} when $\beta=2$.
\subsection*{Acknowledgements}
The author wishes to acknowledge helpful conversations about this work
with J.~P.~Keating and D.~Savin.
\def\Dbar{\leavevmode\lower.6ex\hbox to 0pt{\hskip-.23ex \accent"16\hss}D}
  \def\cprime{$'$} \def\rmi{{\mathrm i}}

\end{document}